\documentclass{amsart}
\usepackage{amssymb,amsfonts,amsmath}
\usepackage{graphicx,xcolor}

\newtheorem{theorem}{Theorem}[section]

\newtheorem{corollary}[theorem]{Corollary}
\newtheorem{proposition}[theorem]{Proposition}

\theoremstyle{definition}
\newtheorem{definition}[theorem]{Definition}
\newtheorem{example}[theorem]{Example}
\newtheorem{assumption}[theorem]{Assumption}
\newtheorem{remark}[theorem]{Remark}

\numberwithin{equation}{section}

\newcommand{\X}{{\mathbb{X}}}
\newcommand{\Y}{{\mathbb{Y}}}
\newcommand{\ZZ}{{\mathbb{Z}}}
\newcommand{\M}{\mathfrak{M}}
\newcommand{\m}{\mathfrak{m}}
\newcommand{\q}{\mathfrak{q}}
\newcommand{\n}{\mathfrak{n}}

\def\epsilon{{\varepsilon}}

\def\phi{{\varphi}}
\let\Psi=\varPsi
\let\Phi=\varPhi
\let\theta=\vartheta
\let\rho=\varrho

\def\LT{\mathop{\rm LT}\nolimits}

\def\LM{\mathop{\rm LM}\nolimits}

\def\HF{\mathop{\rm HF}\nolimits}
\def\HP{\mathop{\rm HP}\nolimits}
\def\ri{\mathop{\rm ri}\nolimits}

\def\Supp{\mathop{\rm Supp}\nolimits}

\def\charac{\mathop{\rm char}\nolimits}

\def\Ann{\mathop{\rm Ann}\nolimits}

\def\Ker{\mathop{\rm Ker}\nolimits}
\def\Im{\mathop{\rm Im}\nolimits}

\def\gr{\mathop{\rm gr}\nolimits}
\def\Tr{\mathop{\rm Tr}\nolimits}
\def\grad{\mathop{\rm grad}\nolimits}

\def\OmegaR{{\Omega^1_{R/K}}}
\def\OmegaRtilde{{\Omega^1_{\widetilde{R}/K}}}
\def\OmegaS{{\Omega^1_{S/K}}}
\def\OmegaRm{{\Omega^m_{R/K}}}
\def\OmegaRtm{{\Omega^m_{\widetilde{R}/K}}}

\let\To=\longrightarrow
\def\TTo#1{\mathop{\longrightarrow}\limits ^{#1}}

\def\deh{^{\rm deh}}

\def\tfrac #1#2{{\textstyle\frac{#1}{#2}}}
\def\tsum{\textstyle\sum\limits}

\def\cocoa{\mbox{\rm
  C\kern-.13em o\kern-.07 em C\kern-.13em o\kern-.15em A}}
\def\apcocoa{\mbox{\rm
A\kern-0.13em p\kern -0.07em C\kern-.13em o\kern-.07 em C\kern-.13em
o\kern-.15em A}}

\begin{document}

\title{Differential Theory of Zero-dimensional Schemes}

%    author one information
\author{Martin Kreuzer}
\address{Faculty of Informatics and Mathematics, University of Passau,
D-94030 Passau, Germany}
\email{martin.kreuzer@uni-passau.de}

%    author two information
\author{Tran N.~K.~Linh}
\address{Department of Mathematics, University of Education - Hue University,
	34 Le Loi, Hue, Vietnam}
\email{tnkhanhlinh@hueuni.edu.vn}

%    author three information
\author{Le N.~Long}
\address{Department of Mathematics, University of Education - Hue University,
	34 Le Loi, Hue, Vietnam}
\curraddr{Faculty of Informatics and Mathematics, University of Passau,
D-94030 Passau, Germany}
\email{lengoclong@dhsphue.edu.vn}

%    \subjclass is required.

\date{\today}

\dedicatory{This paper is dedicated to the memory of Ernst Kunz (1933-2021).}

\begin{abstract}
For a 0-dimensional scheme $\X$ in $\mathbb{P}^n$ over a perfect field~$K$,
we first embed the homogeneous coordinate ring~$R$ into its truncated integral 
closure~$\widetilde{R}$. Then we use the corresponding map from the module of
K\"ahler differentials $\Omega^1_{R/K}$ to $\Omega^1_{\widetilde{R}/K}$ to find
a formula for the Hilbert polynomial $\HP(\Omega^1_{R/K})$ and a sharp bound
for the regularity index $\ri(\Omega^1_{R/K})$. Additionally, we extend this to formulas
for the Hilbert polynomials $\HP(\Omega^m_{R/K})$ and bounds for the regularity indices
of the higher modules of K\"ahler differentials. Next we derive a new 
characterization of a weakly curvilinear scheme~$\X$ which can be checked
without computing a primary decomposition of its vanishing ideal~$I_{\X}$.
Moreover, we prove precise formulas for the Hilbert polynomial of $\Omega^m_{R/K}$
of a fat point scheme~$\X$, extending and settling previous partial results and conjectures.
Finally, we characterize uniformity conditions on~$\X$ using the Hilbert functions of the
K\"ahler differential modules of~$\X$ and its subschemes.
\end{abstract}

\maketitle

%%%%%%%%%%%%%%%%%%%%%%%%%%%%%%%%%%%%%%%%%%%%%%%%%%%%
%
%  Section 1: Introduction
%
%%%%%%%%%%%%%%%%%%%%%%%%%%%%%%%%%%%%%%%%%%%%%%%%%%%%

\section{Introduction}

The study of 0-dimensional subschemes~$\X$ of a projective space $\mathbb{P}^n$, 
in particular of finite set of points, has a long and rich history.
Traditional tools used in this area are the homogeneous vanishing 
ideal~$I_{\X}$ of~$\X$ in $P=K[x_0,\dots,x_n]$, the homogeneous coordinate ring
$R=P/I_{\X}$, its Hilbert function $\HF_{\X}(i)=\dim_K(R_i)$, and the canonical
module~$\omega_R$ of~$R$. In several previous works, the authors introduced and
started using the K\"ahler differential module $\Omega^1_{R/K}$ of~$R$ and even the
entire K\"ahler differential algebra $\Omega^{\bullet}_{R/K} = \Lambda_R(\OmegaR)$
in order to advance this topic (see~\cite{DK1999}, \cite{KLL2019}, and~\cite{KLL2021}).

In the present paper we systematize this approach, extend it by several 
new techniques, solve and unify previous partial results and conjectures,
and derive a number of useful applications to some computational tasks.
Let us describe the individual contributions in more detail.

In Section~2 we recall some basic objects related to a 0-dimensional scheme~$\X$
in the projective space $\mathbb{P}^n$ over a perfect field~$K$, in particular
its homogeneous vanishing ideal $I_\X \subseteq P=K[X_0,\dots,X_n]$, 
its homogeneous coordinate ring $R=P/I_\X$, and its Hilbert function $\HF_\X(i)=
\dim_K(R_i)$ for $i\in\mathbb{Z}$. We always assume that the support of~$\X$
is contained in the affine space $D_+(X_0)$ and denote its affine coordinate
ring by $S=K[X_1,\dots,X_n]/I_\X\deh$. After identifying the homogeneous 
ring of quotients of~$R$ with $S[x_0,x_0^{-1}]$, we introduce a new tool, namely
the embedding of~$R$ into its {\it truncated integral closure} $\widetilde{R}
\cong S[x_0]$.

These techniques come to fruition in Section~3, where we introduce 
and study the K\"ahler differential module $\OmegaR$ and its Hilbert function
$\HF_\OmegaR$. New results here are the description of $\HP(\OmegaR)$
and a sharp bound $\ri(\OmegaR) \le 2\, r_\X + 1$ for the (Hilbert) regularity index
of~$\OmegaR$. Here $r_\X$ is the regularity index of~$\X$, i.e., the first
degree from where on the Hilbert function agrees with the value of the 
Hilbert polynomial (see Prop.~\ref{prop:HKDprops}).
After finding the Hilbert function and the Hilbert polynomial of~$\tilde{R}$
(see Prop.~\ref{prop:OmegaRtilde}), we construct the canonical map
$\Phi:\, \OmegaR \rightarrow \OmegaRtilde$ explicitly, show that $\Ker(\Phi)$
is exactly the torsion submodule of $\OmegaR$, and use it to prove our first main
result, namely the formula
$$
\HP (\OmegaR) \;=\; \deg(\X) + \dim_K(\OmegaS) 
$$
for the Hilbert polynomial of~$\OmegaR$.

Further tools are provided in Section~4 in order to study the torsion
submodule $T\OmegaR$ and other important submodules of~$\OmegaR$.
The Euler derivation $\delta:\, R \rightarrow R$ given by $\delta(f)=
i\, f$ for $f\in R_i$ gives rise to the {\it Euler form}
$\epsilon:\, \OmegaR \rightarrow R$ which satisfies $\epsilon(dx_i)=x_i$
for $i=0,\dots,n$. The Koszul complex over~$\epsilon$ is called the 
{\it Euler-Koszul complex} of~$R$, and the image of $\epsilon^{(2)}:\,
\Lambda^2 \OmegaR \rightarrow \OmegaR$ is called the {\it Koszul submodule}
$U_{R/K} = \langle x_i dx_j - x_j dx_i \mid 0\le i<j\le n\rangle$ of~$\OmegaR$.
Both $T\OmegaR$ and $U_{R/K}$ are contained in $\Ker(\epsilon)$
(see Props.~\ref{prop:EulerProps} and~\ref{prop:KoszulEuler}). The Koszul
submodule is equal to $\Ker(\epsilon)$ if $\charac(K)$ is zero or 
$\ge 2 r_\X+1$, and in the other cases the (small) difference between
the two submodules is laid out in detail in Prop.~\ref{prop:KerEuler}.

Next up, we apply the preceding results in Section~5 to the higher
K\"ahler differential modules $\OmegaRm$, where $m\ge 1$. After recalling
the presentation and computation of~$\OmegaRm$ 
(see Prop.~\ref{prop:PresentationOmegaRm}), we use the embedding
$R\hookrightarrow \tilde{R}$ to prove the formula
$$
\HP(\OmegaRm) \;=\; \dim_K(\Omega^m_{S/K}) + \dim_K( \Omega^{m-1}_{S/K})
$$
for the Hilbert polynomial of~$\OmegaRm$ and the sharp bound
$\ri(\OmegaRm) \le r_\X+m$ for its regularity index 
(see Prop.~\ref{prop:HigherKDM}).

In final three sections we apply and extend these algebraic results
to characterize and compute geometric properties of the scheme~$\X$.
In Section~6 we look at (weakly) curvilinear schemes which are defined by the
property that the maximal ideals of their local rings are unigenerated.
Notice that this generalizes slightly the condition that these local
rings are of the form $K[z]/\langle p^k\rangle$ with an irreducible
polynomial~$p$. Special cases are, of course, smooth schemes for which it is
known that they can be characterized by $\Omega^1_{S/K} = 0$, or, equivalently, by
$\HP(\OmegaR) = \deg(\X)$ (see Prop.~\ref{prop:SmoothnessCriterion}). 
Note that this allows us to check smoothness of~$\X$ without
computing the primary decomposition of~$I_\X$.
In a similar vein, we characterize non-smooth weakly curvilinear schemes by
$\Omega^1_{S/K} \ne 0$ and $\Omega^m_{S/K}=0$ for $m\ge 2$. Equivalently, the
scheme~$\X$ is weakly curvilinear, but not smooth if and only if
$$
\HP(\OmegaR) > \deg(\X),\; \HP(\Omega^2_{R/K})= \HP(\OmegaR) -\deg(\X),
\; \HP(\Omega^m_{R/K})=0\; \hbox{\rm for}\, m\ge 3
$$
(see Prop.~\ref{prop:CharCurvilinearKD}). Thus we can check algorithmically
whether~$\X$ is weakly curvilinear without computing a primary decomposition.
Using the dimensions of the residue 
class fields of the local rings of~$\X$, we can also write down explicit
formulas for the Hilbert polynomials of~$\OmegaR$ and $\Omega^2_{R/K}$ in 
the weakly curvilinear case (see Cor.~\ref{cor:ExplicitCurvilinear}).

The next application is contained in Section~7 where we consider fat 
point schemes~$\X$. They are defined by vanishing ideals of the form
$I_\X = I_{p_1}^{m_1} \cap\cdots\cap I_{p_t}^{m_t}$ where
the points $p_i$ in the support of~$\X$ are assumed to be $K$-rational
and $m_i\ge 1$. The case of a reduced scheme, i.e., the case
$m_1=\cdots =m_t=1$, is easily characterized by $\HP(\OmegaR)=t$
and $\HP(\OmegaRm)=0$ for $m\ge 2$ (see Remark~\ref{rem:FatPointProps}).
The general case reduces to the study of K\"ahler differential modules
of rings of the form $S=A/\q^k$, where $A=K[X_1,\dots,X_n]$ and 
$\q = \langle X_1,\dots,X_n\rangle$. Here we obtain explicit formulas
for the values of the Hilbert function of $\Omega^m_{S/K}$ which only
depend on the value of $\delta = \dim_K(d\q^k \wedge \Omega^{m-1}_{A/K})$
(see Prop.~\ref{prop:K-HFOmegaSm}). 

If $\charac(K)=0$ or $\charac(K)>k$, we succeed in determining this
value~$\delta$ using a number of subtle and detailed arguments
(see Props.~\ref{prop:LTgradti}, \ref{prop:GensviJ}, and~\ref{prop:LTviJ}).
Thus we are able to derive the explicit formula
$$
\dim_K(\Omega^m_{S/K}) \;=\; \tbinom{n}{m}\tbinom{n+k-2}{n} +
\tbinom{m+k-2}{m}\tbinom{n+k-2}{n-m-1}
$$
for $S=K[X_1,\dots,X_n]/\langle X_1,\dots,X_n\rangle^k$ and $\charac(K)=0$
or $\charac(K)>k$ (see Thm.~\ref{thm:K-DimOmegaSm}).
Under a more stringent assumption about the characteristic of~$K$, we 
also provide another proof based on the exactness of the Euler-Koszul complex
(see Prop.~\ref{prop:KoszulExactSeq} and Remark~\ref{rem:AltProof}).

Finally, we combine everything and get explicit formulas for the
Hilbert polynomials of the modules $\OmegaRm$ if $\charac(K)=0$
or $\charac(K)> \max\{m_1,\dots,m_t\}$. They generalize and vastly improve
the partial results in~\cite{KLL2019} and~\cite{KLL2021}.

To cap the paper off, we look at differential characterizations
of uniformity properties of~$\X$ in Section~8. One of them, the 
Cayley-Bacharach property of degree~$d$, requires that every hypersurface
of degree~$d$ which contains a subscheme~$\Y$ of~$\X$ of colength one
automatically contains~$\X$. It can be characterized by the degrees
of the minimal separators of~$\X$ (see Remark~\ref{rem:FirstCBP}).
Our key result here is that the vanishing ideal $I_{\Y/\X}$ of a 
subscheme~$\Y$ in~$\X$ satisfies $(dI_{\Y/\X})_{\alpha_{\Y/\X}} \ne 0$
in its initial degree $\alpha_{\Y/\X}$ provided that $\charac(K)=0$
or $\charac(K)> r_\X$ (see Prop.~\ref{prop:OmegaSubsch}).
This allows us to characterize schemes having the Cayley-Bacharach property
of degree~$d$ through the Hilbert function of their K\"ahler differential 
modules (see Cor.~\ref{cor:CBP-Omega}). More generally, we consider
higher uniformity properties of~$\X$, called $(i,j)$-uniformity, and
characterize them through the Hilbert function of $\OmegaR$ as well 
(see Cor.~\ref{cor:Uniformity-Omega}).

Many examples in this paper illustrate our results. They were computed
using the computer algebra system ApCoCoA (see~\cite{ApCoCoA}).
For the notation and basic definitions we use, we refer to~\cite{KR2000},
\cite{KR2005}, and~\cite{KR2016}. Our main reference for results about 
K\"ahler differential modules is the fundamental work~\cite{Kun1986} of our late
teacher, mentor, and supporter {\it Ernst Kunz} (1933-2021) to whom we dedicate
this paper with deep gratitude.

%%%%%%%%%%%%%%%%%%%%%%%%%%%%%%%%%%%%%%%%%%%%%%%%%%%%
%
%  Section 2: Zero-Dimensional Schemes
%
%%%%%%%%%%%%%%%%%%%%%%%%%%%%%%%%%%%%%%%%%%%%%%%%%%%%

\section{Zero-Dimensional Schemes}\label{sec:0dimSchemes}

In the following we let~$\X$ be a 0-dimensional scheme in projective $n$-space
$\mathbb{P}^n$ over a perfect field~$K$. The homogeneous coordinate ring of~$\mathbb{P}^n$
is the polynomial ring $P=K[X_0,\dots,X_n]$, and the homogeneous coordinate ring of~$\X$
is $R = P/I_\X$, where $I_\X$ denotes the homogeneous vanishing ideal of~$\X$.
Recall that~$I_\X$ is a saturated homogeneous ideal, and that $R$ is a standard graded
1-dimensional Cohen-Macaulay ring. In particular, there exists a homogeneous non-zerodivisor
in~$R$. To simplify the presentation, we make the following assumption.

\begin{assumption}
Suppose that no point of the support of~$\X$ is contained in the hyperplane $\mathcal{Z}(X_0)$.
\end{assumption}

Usually, if this assumption does not hold, it suffices to perform a homogeneous linear
change of coordinates. This will not affect any results in this paper. In the case of a finite
field~$K$, it may be necessary to perform a small base field extension in order to find
a suitable homogeneous linear change of coordinates. This will not affect any results either.
As a consequence of this assumption, we get that~$\X$ is contained in the affine space
$\mathbb{A}^n \cong D_+(X_0)$, and that the image $x_0$ of~$X_0$ in~$R$ is a non-zerodivisor.

The affine coordinate ring of~$\X$, viewed as a subscheme of~$\mathbb{A}^n$, is given by
$S = R/ \langle x_0- 1\rangle \cong K[X_1,\dots,X_n] / I_\X\deh$, where $I_\X\deh$
is the dehomogenization
of~$I_\X$ with respect to~$X_0$. The ring~$S$ is a 0-dimensional affine $K$-algebra and hence
a finite dimensional $K$-vector space. The following functions are central to the study
of 0-dimensional subschemes of projective spaces.

\begin{definition}
In the above setting, the map $\HF_\X:\; \ZZ \To \ZZ$ given by $\HF_\X(i)= \dim_K(R_i)$
for all $i\in\ZZ$ is called the {\bf Hilbert function} of~$\X$.

Its first difference function $\Delta\HF_\X:\; \ZZ \To\ZZ$ which is given by
$\Delta\HF_\X(i) = \HF_\X(i) - \HF_\X(i-1)$ for $i\in\ZZ$ is called the {\bf Castelnuovo
function} of~$\X$.
\end{definition}

The Hilbert function contains a number of useful invariants of the embedding of~$\X$
into~$\mathbb{P}^n$. Let us collect some of its basic properties.

\begin{remark}
The Hilbert function of a 0-dimensional scheme~$\X$ in $\mathbb{P}^n$ satisfies
$$
1 = \HF_\X(0) < \HF_\X(1) < \cdots < \HF_\X(r_\X) = \HF_\X(r_\X + 1) = \cdots = \deg(\X)
$$
where $r_\X$ is called the (Hilbert) {\bf regularity index} of~$\X$ and $\deg(\X)$
is the degree of~$\X$, i.e., the value of the (constant) Hilbert polynomial of~$R$.
\end{remark}

The structure of the ring~$R$ in degrees $\ge r_\X$ can be described as follows.

\begin{proposition}\label{prop:RinLargeDegrees}
For $i\ge r_\X$, the map $\epsilon_i:\; R_i \To S$ given by $f\mapsto f\deh$ is
an isomorphism of $K$-vector spaces.
\end{proposition}

\begin{proof}
Since~$x_0$ is a non-zerodivisor of~$R$, the definition of~$r_\X$ shows that
the multiplication map $\mu_{x_0}:\; R_i \To R_{i+1}$ is an isomorphism
for $i\ge r_\X$. Together with the observation that $\dim_K(R_i) = \dim_K(S) = \deg(\X)$
for $i\ge r_\X$, it follows that it suffices to prove the claim for $i\gg 0$.
Let $\{b_1,\dots,b_r\}$ be a $K$-vector space basis of~$S$.
Using the fact that the map $\delta:\; R \To R/ \langle x_0-1\rangle \cong S$ is surjective
and given by dehomogenization, we can find for every $j\in\{1,\dots,r\}$
a homogeneous element $f_j\in R$ such that $\delta(f_j)=f_j\deh = b_j$. By multiplying
each $f_j$ with an appropriate power of~$x_0$, we may assume that all elements $f_j$
have the same degree~$i$. Thus the map $\epsilon_i$ is surjective, hence bijective,
and the claim follows.
\end{proof}

Subsequently, we will also use the following ring. The graded $K$-algebra
$$
Q^h(R) \;=\; \{ \tfrac{a}{b} \mid a,b\in R; \; b \hbox{\ \rm homogeneous non-zerodivisor }\}
$$
is called the (full) {\bf homogeneous ring of quotients} of~$R$. We can simplify its
description as follows.

\begin{proposition}\label{prop:QhofR}
Let~$R$ be the homogeneous coordinate ring of a 0-dimensional
scheme~$\X$ in~$\mathbb{P}^n$ as above.
\begin{enumerate}
\item[(a)] The canonical inclusion $R_{x_0} \rightarrow Q^h(R)$ is bijective.

\item[(b)] The map $\phi:\; Q^h(R) \To S[x_0,x_0^{-1}]$ given by $\phi(f/x_0^i)
= f\deh\, x_0^{d-i}$ for $f\in R_d$ and $i\ge 0$ is an isomorphism of graded $K$-algebras.
\end{enumerate}
\end{proposition}

\begin{proof}
To prove (a), we let $a\in R_d$ and $b\in R_e$ be homogeneous elements of~$R$ 
of degrees $d,e\in\mathbb{N}$, and assume that~$b$ is a non-zerodivisor.
We need to show that there exist $i\ge 0$ and a homogeneous element $f\in R_i$ such that
$a/b = f/x_0^i$. First we choose~$i$ such that $i\ge r_\X + e -d$.
Note that the multiplication map $\mu_b:\, R_k \To R_{k+e}$ is bijective for every $k\ge r_\X$,
because~$b$ is a non-zerodivisor. Hence $d+i\ge r_\X+e$ implies that
the element $a\, x_0^i \in R_{d+i}$ is of the form $a\, x_0^i = b\, f$ for some
$f\in R_j$, where $j=d+i-e$. Thus we get $a/b = f/x_0^i$, as desired.

It remains to prove~(b). Given $f\in R_d$ and $i\ge 0$, we let $j=r_\X +i$. Then we
have $x_0^j\, (f/x_0^i) \in R_{r_\X +d}$ and $x_0^j\, \phi(f/x_0^i) = \phi(x_0^j\, (f/x_0^i)) =
f\deh x_0^{r_\X+d}$. Together with Proposition~\ref{prop:RinLargeDegrees}, it follows
that $\phi_{d-i}:\; Q^h(R)_{d-i} \To S x_0^{d-i}$ is an isomorphism of $K$-vector spaces.
As~$\phi$ is clearly a homomorphism of graded rings, the claim follows.
\end{proof}

Using the identifications $Q^h(R)\cong R_{x_0} \cong S[x_0,x_0^{-1}]$ provided by
this proposition, we now examine the integral closure of~$R$ in~$Q^h(R)$.

\begin{proposition}\label{prop:RandRtilde}
Let~$R$ be the homogeneous coordinate ring of a 0-dimensional
scheme~$\X$ in~$\mathbb{P}^n$ as above.
\begin{enumerate}
\item[(a)] The map $\psi:\; R \To S[x_0]$ given by
$\psi(f) = f\deh\, x_0^d$ for $f\in R_d$ is an injective graded $K[x_0]$-algebra homomorphism.
\item[(b)] For $i\ge r_\X$, the map $\psi_i:\; R_i \To Sx_0^i$ is an isomorphism
of $K$-vector spaces.

\item[(c)] The subring~$\widetilde{R}$ of~$Q^h(R)$ generated by all homogeneous
elements of non-negative degree satisfies $\widetilde{R}\cong S[x_0]$
and is contained in the integral closure of~$R$ in~$Q^h(R)$.

\item[(d)] The ring $\widetilde{R}$ is the integral closure of~$R$ 
in~$Q^h(R)$ if and only if $S$ is a reduced ring.
\end{enumerate}
\end{proposition}

\begin{proof}
Claims (a) and~(b) follow immediately from Propositions~\ref{prop:RinLargeDegrees}
and~\ref{prop:QhofR}. To prove~(c), we let $f\in Q^h(R)$ be a homogeneous element 
of degree $d\ge 0$. Using the isomorphism $\phi:\; Q^h(R) \To S[x_0,x_0^{-1}]$ 
given in Proposition~\ref{prop:QhofR},
we write $\phi(f) = x_0^d\, g$ with $g\in S$. By definition, the map~$\phi$ 
identifies~$\widetilde{R}$ with $S[x_0]$.
Since~$S$ is a finite $K$-algebra, there exists
a relation $g^m + a_{m-1} g^{m-1} + \cdots + a_1 g_1 + a_0 = 0$ with $a_i \in K$. 
Then we apply~$\phi^{-1}$ to both sides of 
$$
x_0^{dm} \, g^m + a_{m-1} x_0^d (x_0^{(m-1)d} g^{m-1}) + \cdots 
+ a_1 x_0^{(m-1)d} ( x_0^d g )  + a_0 x_0^{dm}  \;=\; 0
$$
and get $f^m + a_{m-1} x_0^d f^{m-1} + \cdots + a_1 x_0^{(m-1)d} f + a_0 x_0^{dm} =0$.
Thus $a_i x_0^j \in R$ implies that~$f$ is integral over~$R$.

Finally, we prove~(d). If~$S$ is a reduced ring, then $\widetilde{R}$
is the integral closure of~$R$ in~$Q^h(R)$ by
\cite{Bor1989}, Ch.~V, \S 2, Prop.~9 and \S 8, Prop.~20.
Conversely, if~$S$ contains a non-zero nilpotent element $f\in S$
with $f^k=0$ for some $k\ge 1$, then we consider the homogeneous element 
$g = \phi^{-1}(fx_0^{-1}) \in Q^h(R)_{-1}$. By definition, we have $g\notin \widetilde{R}$. 
Moreover, from $\phi(g^k) = f^k x_0^{-k} =0$, it follows that $g^k=0$, and therefore~$g$ is integral 
over~$R$.
\end{proof}

The ring~$\widetilde{R}$ will prove very useful and deserves a name.

\begin{definition}
Let~$R$ be the homogeneous coordinate ring of a 0-dimensional
scheme~$\X$ in~$\mathbb{P}^n$ as above. Then the subring~$\widetilde{R}$ of~$Q^h(R)$ generated by
all homogeneous elements of non-negative degree is called the {\bf truncated integral closure}
of~$R$.
\end{definition}

Frequently, we will use the identification $\widetilde{R} \cong S[x_0]$ 
given by $f \mapsto f^{\rm deh}\, x_0^d$ for a homogeneous element $f\in \widetilde{R}_d$
of degree $d\ge 0$.

\bigbreak
%%%%%%%%%%%%%%%%%%%%%%%%%%%%%%%%%%%%%%%%%%%%%%%%%%%%
%
%  Section 3: The Kaehler Differential Module
%
%%%%%%%%%%%%%%%%%%%%%%%%%%%%%%%%%%%%%%%%%%%%%%%%%%%%

\section{The K\"ahler Differential Module}\label{sec:KDmodule}

In the following we continue to use the assumptions and notation
introduced above. In particular, let $K$ be a perfect field, and let
$\X$ be a 0-dimensional subscheme of~$\mathbb{P}^n$ with homogeneous
coordinate ring $R = P/I_\X$ for $P=K[X_0,\dots,X_n]$. For $i=0,\dots,n$,
we denote the image of~$X_i$ in~$R$ by~$x_i$, and we assume that~$x_0$
is a non-zerodivisor in~$R$. Furthermore, we let $\mathfrak{m}=
\langle x_0,\dots,x_n\rangle$ be the homogeneous maximal ideal of~$R$.
As the following $R$-module is the main object of study in this section,
we recall its (well-known) definition.

\begin{definition}\label{def:KDdef}
Let $J$ be the kernel of the multiplication map $\mu: R \otimes_K R \To R$.
Then $\OmegaR = J/J^2$ is an $R$-module via
$r\cdot \sum_i a_i \otimes b_i = \sum_i ra_ib_i$ for $r,a_i,b_i \in R$.
It is called the {\bf module of K\"ahler differentials} of $R/K$, or the
{\bf K\"ahler differential module} of~$R/K$.

The map $d_{R/K}:\; R \To \OmegaR$ defined by $d_{R/K}(f) = f\otimes 1 - 1 \otimes f + J^2$
is called the {\bf universal derivation} of~$R/K$.
\end{definition}

Notice that the map $d_{R/K}$ is indeed a derivation of $R/K$, i.e.,
it is $K$-linear and satisfies the product rule. If the algebra $R/K$ is clear from the
context, we will usually simply write $df$ instead of $d_{R/K}(f)$.
Let us recall some of the basic properties of the K\"ahler differential module
of~$R/K$.

\begin{proposition}\label{prop:presentOmega1}
Let $R$ be the homogeneous coordinate ring of a 0-dimensional scheme~$\X$ in~$\mathbb{P}^n$
as above.
\begin{enumerate}
\item[(a)] The $R$-module $\OmegaR$ is positively graded and generated by the homogeneous
elements $\{dx_0,\dots, dx_n\}$ of degree~1.

\item[(b)] There is a homogeneous short exact sequence of graded $R$-modules
$$
I_\X / I_\X^2  \;\TTo{\alpha}\; R^{n+1}(-1)  \;\TTo{\beta}\;  \OmegaR \;\To\; 0
$$
where~$\alpha$ is given by $\alpha(f+I_\X^2) = (\tfrac{\partial f}{\partial x_0}, \dots,
\tfrac{\partial f}{\partial x_n})$ and where~$\beta$ is given by $\beta(e_i)=dx_i$
for $i=0,\dots,n$.

\item[(c)] The graded $R$-module $\OmegaR$ has a presentation
$$
\OmegaR \;\cong\;  \Omega^1_{P/K} / (I_\X \cdot \Omega^1_{P/K} + dI_\X)
$$
where $\Omega^1_{P/K} = PdX_0 \oplus \cdots \oplus PdX_n \cong P^{n+1}(-1)$
is a graded free $P$-module and $dI_\X = \langle \tfrac{\partial f}{\partial X_0} dX_0
+ \cdots + \tfrac{\partial f}{\partial X_n} dX_n \mid f\in I_\X \rangle$.
\end{enumerate}
\end{proposition}

\begin{proof}
Claims (a) and (b) follow from \cite{Kun1986}, Propositions~4.12 and~4.17
and claim (c) follows from ibid., Proposition~4.19.
\end{proof}

Every finitely generated graded $R$-module has a Hilbert function, a constant
Hilbert polynomial, and a regularity index. This yields the following definition.

\begin{definition}
In the above setting, let $\OmegaR$ be the K\"ahler differential module of~$R$.
\begin{enumerate}
\item[(a)] The map $\HF_\OmegaR:\; \ZZ \To \ZZ$ given by $\HF_\OmegaR(i) =
\dim_K(\OmegaR)_i$ for $i\in \ZZ$ is called the {\bf Hilbert function} of~$\OmegaR$.

\item[(b)] The number $\HP(\OmegaR) = \HF_\OmegaR(i)$ for $i \gg 0$ is called the
{\bf Hilbert polynomial} of~$\OmegaR$.

\item[(c)] The number $\ri(\OmegaR) = \min\{ i\in\ZZ \mid \HF_\OmegaR(j)=
\HP(\OmegaR) \hbox{\ \rm for \ }j\ge i\}$ is called the {\bf regularity index} of~$\OmegaR$.
\end{enumerate}
\end{definition}

The following proposition shows that the Hilbert polynomial of~$\OmegaR$
is well-defined and that there is a bound for its regularity index.

\begin{proposition}\label{prop:HKDprops}
Let $R$ be the homogeneous coordinate ring of a 0-dimensional scheme~$\X$ 
in~$\mathbb{P}^n$ as above.
\begin{enumerate}
\item[(a)] For $i\le 0$, we have $\HF_\OmegaR(i)=0$.

\item[(b)] For $i\ge r_\X +1$, the multiplication map
$\mu_{x_0}:\; (\OmegaR)_i \To (\OmegaR)_{i+1}$ is surjective.
In particular, we have $\ri(\OmegaR)\ge r_\X+1$.

\item[(c)] If $\ri(\OmegaR) > r_\X+1$ then we have
$$
  \HF_\OmegaR (r_\X + 1) > \cdots >  \HF_\OmegaR (\ri(\OmegaR))
$$

\item[(d)] For $i\ge 2r_\X+1$, the multiplication map
$\mu_{x_0}:\; (\OmegaR)_i \To (\OmegaR)_{i+1}$ is bijective.
In particular, we have $\ri(\OmegaR) \le 2 r_\X +1$.
\end{enumerate}
\end{proposition}

\begin{proof}
Claim (a) follows from the fact that $\OmegaR$ is generated by homogeneous
elements of degree~1. To prove (b), we let $i\ge r_\X+1$ and $w\in (\OmegaR)_{i+1}$.
Then we can write $w= f_0 dx_0 + \cdots + f_n dx_n$ with $f_0,\dots,f_n \in R_i$.
Since $i\ge r_\X+1$, Proposition~\ref{prop:RinLargeDegrees} implies that there
exist elements $g_j\in R_{i-1}$ such that $f_j = x_0 g_j$ for $j=0,\dots,n$, and
hence $w = x_0\, (g_0 dx_0 + \cdots +g_n dx_n)$ is in the image of~$\mu_{x_0}$.

Claim~(c) is a consequence of~(b) and the definition of $\ri(\OmegaR)$.
It remains to show~(d). Let $w = f_0 dx_0 + \cdots + f_n dx_n$ be an element
in the kernel of~$\mu_{x_0}$, where $f_j \in R_{i-1}$ for $j=0,\dots,n$.
By Proposition~\ref{prop:presentOmega1}, we have an exact sequence of graded $R$-modules
$$
0 \;\To\; \mathcal{G}(-1) \;\To\;  R^{n+1}(-1) \;\TTo{\phi}\; \OmegaR \;\To\; 0
$$
where $\mathcal{G}$ is generated by the tuples $(\frac{\partial f}{\partial x_0}, \dots,
\frac{\partial f}{\partial x_n})$ such that $f\in I_\X$ and where
$\phi(f_0,...,f_n) = f_0 dx_0 + \cdots + f_n dx_n$ for all $f_0,...,f_n\in R$.
Since $I_\X$ is generated in degrees $\le r_\X+1$ by~\cite{GM1984}, Proposition~1.1,
the graded $R$-module $\mathcal{G}$ is generated in degrees $\le r_\X$.
So, let $\{v_1,...,v_k\}$ be a homogeneous system of generators of~$\mathcal{G}$,
where $\deg(v_j)\le r_\X$. Given an element $w = f_0 dx_0 + \cdots + f_n dx_n
\in \Ker(\mu_{x_0})$ with $f_j\in R_{i-1}$ for $j=0,\dots,n$, we have
$$
(x_0 f_0, \dots, x_0 f_n) \;\in\; \Ker(\phi)  \;=\;
[\mathcal{G}(-1)]_{i+1}  \;=\;  \mathcal{G}_i
$$
Thus we may write $(x_0f_0,...,x_0f_n) = g_1 v_1 + \cdots + g_k v_k$
with homogeneous elements $g_j\in R$ of degree $\deg(g_j) = i - \deg(v_j)
\ge i - r_\X \ge r_\X+1$. Now Proposition~\ref{prop:RinLargeDegrees} allows us
to write $g_j= x_0 h_j$ with $h_j\in R$ and we obtain
$$
(f_0, \dots, f_n) \;=\; h_1 v_1 + \cdots + h_k v_k \in \mathcal{G}_{i-1}
$$
This yields $w=	f_0 dx_0 + \cdots + f_n dx_n = 0$ in $\OmegaR$, and the proof is complete.
\end{proof}

The graded $R$-module
$$
T\OmegaR \;=\;  \{w\in \OmegaR \mid rw =0\;
\mbox{ for some homogeneous non-zerodivisor}\ r\in R \}
$$
is called the {\bf torsion submodule} of $\OmegaR$.
We can describe it as follows.

\begin{proposition} \label{prop:TorsionMod}
Let $R$ be the homogeneous coordinate ring of~$\X$ as above.
\begin{enumerate}
\item[(a)] We have $T\OmegaR = \{\, w\in \OmegaR \mid x_0^iw =0$ for some  
$i\ge 1 \,\}$.
	
\item[(b)] We have $\HP(T\OmegaR)=0$ and $\ri(T\OmegaR)\le 2r_\X+1$.
\end{enumerate}
\end{proposition}

\begin{proof}
Let $w\in T\OmegaR$ be a non-zero homogeneous element with
$fw=0$ for a homogeneous non-zerodivisor $f\in R_i$, $i\ge 0$.
Then $f\deh \in S$ is a non-zerodivisor. Since $S$ is an Artinian semilocal ring,
it consists of only units and zero-divisors, and so $f\deh \in S$ is a unit.
Let $u\in S$ be such that $uf\deh =1$ and $g=\psi_{r_\X}^{-1}(ux_0^{r_\X})$,
where $\psi_{r_\X}: R_{r_\X}\rightarrow Sx_0^{r_\X}$
is the isomorphism of $K$-vector spaces in Proposition~\ref{prop:RandRtilde}.b.
We have $gf = x_0^{r_\X+i}$, and so
$x_0^{r_\X+i}w = gfw =0$.
This proves claim (a).
Claim (b) follows from (a) and Proposition~\ref{prop:HKDprops}.d.
\end{proof}

Let us show by the following example that the above upper
bound for the regularity index of $T\OmegaR$ is a sharp bound.

\begin{example}
Let $\X$ and~$\Y$ be two sets of five $K$-rational points
in $\mathbb{P}^2$ such that~$\X$ consists of three points on a line
and two points on another line (two lines can intersect at one point of $\X$), 
and such that~$\Y$ is contained in a non-singular conic.
By~$R_\Y$ we denote the homogeneous coordinate of~$\Y$.  
Then the Hilbert functions of~$R$ and $R_\Y$ agree, namely
$
\HF_\X = \HF_\Y: \ 1 \ 3\ 5\ 5\cdots,
$
and so $r_\X = r_\Y=2$. However, the Hilbert functions
of $\OmegaR$ and $\Omega^1_{R_\Y/K}$
as well as of $T\OmegaR$ and $T\Omega^1_{R_\Y/K}$
are different:
$$
\begin{aligned}
\HF_{\OmegaR}&:\ 0\ 3\ 8\ 10\ \textbf{7}\ 5\ 5\cdots, &\quad&
\HF_{\Omega^1_{R_\Y/K}}\ \,:\ 0\ 3\ 8\ 10\ \textbf{6}\ 5\ 5\cdots,\\
\HF_{T\OmegaR}&:\ 0\ 0\ 3\ 5\ \textbf{2}\ 0\ 0\cdots, &\quad&
\HF_{T\Omega^1_{R_\Y/K}}:\ 0\ 0\ 3\ 5\ \textbf{1}\ 0\ 0\cdots.
\end{aligned}
$$
In particular, we get $\ri(\OmegaR) = \ri(T\OmegaR) =
\ri(\Omega^1_{R_\Y/K}) = \ri(T\Omega^1_{R_\Y/K}) = 5 = 2r_\X+1$, 
and hence the upper bounds for the regularity index
of $\OmegaR$ and $T\OmegaR$ given in Propositions~\ref{prop:HKDprops}
and~\ref{prop:TorsionMod} are sharp.
Moreover, these upper bounds are also sharp for the non-reduced 
0-dimensional scheme $\Y'\subseteq \mathbb{P}^2$
whose support comprises five points $(1:0:1)$, $(1:1:2)$, $(1:2:2)$, 
$(1:3:1)$, and $(1:1:0)$ on a non-singular conic, 
in which only the last point $(1:1:0)$ is non-reduced with its ideal
$I_{(1:1:0)} = \langle X_1-X_0, X_2^2\rangle \subseteq K[X_0,X_1,X_2]$
and ${\rm char}(K)\ne 2,3$, because a calculation gives
$\HF_{\Y'}:\ 1\ 3\ 6\ 6\cdots$, $r_{\Y'} =2$, and
$$
\begin{aligned}
	&\HF_{\Omega^1_{R_{\Y'}/K}}:\ 0\ 3\ 9\ 14\ 9\ 7\ 7\ \cdots,\
	\ri(\Omega^1_{R_{\Y'}/K})=5 = 2r_{\Y'}+1,\\
	&\HF_{T\Omega^1_{R_{\Y'}/K}}:\ 0\ 0\ 2\ 7\ 2\ 0\ 0\ \cdots,\
	\ri(T\Omega^1_{R_{\Y'}/K})= 5 = 2r_{\Y'}+1.
\end{aligned}
$$
\end{example}

Next we consider the module of K\"ahler differentials of 
the truncated integral closure $\widetilde{R} \cong S[x_0]$.
Its Hilbert function, Hilbert polynomial, and regularity index 
can be described as follows.

\begin{proposition}\label{prop:OmegaRtilde}
Let $\widetilde{R}\cong S[x_0]$ be the truncated integral closure of~$R$
in its homogeneous ring of quotients.
\begin{enumerate}
\item[(a)] We have an isomorphism of graded $K[x_0]$-modules
$$
\phi:\; \OmegaRtilde \cong S[x_0] dx_0 \oplus K[x_0] \otimes_K\OmegaS.
$$

\item[(b)] The Hilbert function $\HF_\OmegaRtilde:\; \ZZ \To \ZZ$
of $\OmegaRtilde$ is given by $\HF_\OmegaRtilde(i) = 0$ for $i<0$,
by $\HF_{\OmegaRtilde}(0) = \dim_K(\OmegaS)$, and by $\HF_{\OmegaRtilde}(i) =
\deg(\X) + \dim_K(\OmegaS)$ for $i>0$.

\item[(c)] We have $\HP(\OmegaRtilde)= \deg(\X) + \dim_K(\OmegaS)$ and
$\ri(\OmegaRtilde) = 1$.
\end{enumerate}
\end{proposition}

\begin{proof}
Claim~(a) follows from~\cite{Kun1986}, Formula~4.11.a, and from 
$\Omega^1_{K[x_0]/K} \cong K[x_0] dx_0$. 
Part~(b) is an immediate consequence of~(a), and~(c) follows from~(b).
\end{proof}

Next we use the injective homomorphism of graded
$K$-algebras $\phi:\; R \to \widetilde{R}\cong S[x_0]$
given by $f\mapsto f\deh x_0^k$ for $f\in R_k$ to compare the K\"ahler differential
modules of~$R$ and~$\widetilde{R}$.

\begin{proposition}\label{prop:CompareKahlerMod}
Let $R$ be the homogeneous coordinate ring of a 0-dimensional scheme~$\X$
in~$\mathbb{P}^n$ as above, and let $\widetilde{R}\cong S[x_0]$ be the
truncated integral closure of~$R$ in its homogeneous ring of quotients.
\begin{enumerate}
\item[(a)] The homomorphism $\phi:\; R \to \widetilde{R}\cong S[x_0]$
induces a homomorphism of graded $K[x_0]$-modules
$\Phi:\; \OmegaR \To \OmegaRtilde \cong
S[x_0] dx_0 \oplus K[x_0] \otimes_K \OmegaS $ which is given by
$\Phi(f dx_i) = f\deh x_i x_0^k dx_0 + x_0^{k+1} \otimes f\deh dx_i$
for $i=1,\dots,n$, and by $\Phi(f dx_0) = f\deh x_0^k dx_0$ for $f\in R_k$.

\item[(b)] For $i\ge 2r_\X+1$, the map $\Phi_i:\; (\OmegaR)_i \To
(\OmegaRtilde)_i$ is an isomorphism of $K$-vector spaces.

\item[(c)] We have $\HP(\OmegaR) = \deg(\X) + \dim_K(\OmegaS)$
and $\ri(\OmegaR) \le 2 r_\X +1$.

\item[(d)] We have $T\OmegaR = \Ker(\Phi)$.
\end{enumerate}
\end{proposition}

\begin{proof}
To prove~(a), it suffices to use~$\phi$ and calculate 
$\Phi(fdx_i) = \phi(f)d\phi(x_i) = f\deh x_0^k d(x_ix_0) 
= f\deh x_i x_0^k dx_0 + x_0^{k+1} \otimes f\deh dx_i$ for $i=1,\dots,n$
and $\Phi(fdx_0) = f\deh x_0^k dx_0$.

Next we show~(b). Since the multiplication maps 
$\mu_{x_0}:\; (\OmegaR)_i \To (\OmegaR)_{i+1}$
and  $\mu_{x_0}:\; (\OmegaRtilde)_i \To (\OmegaRtilde)_{i+1}$ 
are isomorphisms for $i\ge 2r_\X +1$,
it suffices to prove that $\Psi:=\Phi_{2r_\X+1}:\; (\OmegaR)_{2r_\X+1} \To
(\OmegaRtilde)_{2r_\X+1}$ is an isomorphism. First we prove that~$\Psi$ is injective.
The inclusions $R\subset \widetilde{R} \subset Q^h(R) = R_{x_0}$ induce canonical
homomorphisms of graded $K[x_0]$-algebras
$$
\OmegaR \;\To\; \OmegaRtilde \;\To\; \Omega^1_{Q^h(R)/K} \;\cong\; \Omega^1_{R_{x_0}/K}
\;\cong\; (\OmegaR)_{x_0}
$$
whose composition is the canonical map to the localization of~$\OmegaR$ 
in the element~$x_0$. If $w\in (\Omega_R)_{2r_\X+1}$ is in the kernel 
of~$\Psi$ then its image $\frac{w}{1}$ in $(\OmegaR)_{x_0}$ is zero. 
Hence there exists a number $j\ge 0$ such that $x_0^j\,w = 0$,
and Proposition~\ref{prop:HKDprops}.d yields $w=0$.

It remains to prove that~$\Psi$ is surjective. In view of Proposition~\ref{prop:OmegaRtilde}.a,
we have to consider two cases. For an element of the form $f x_0^{2r_\X} dx_0$, we let
$F\in R_{r_\X}$ be a preimage of $f x_0^{r_\X}$ under the isomorphism $R_{r_\X} \cong S x_0^{r_\X}$
in Proposition~\ref{prop:RandRtilde}.b. Then we have $\Psi(F x_0^{r_\X} dx_0) = F\deh x_0^{2r_\X} dx_0
= f x_0^{2r_\X} dx_0$. For an element of the form $x_0^{2r_\X+1} \otimes fdg$ with $f,g\in S$, we let
$F,G\in R_{r_\X}$ be preimages of $fx_0^{r_\X}$ and $gx_0^{r_\X}$ under the isomorphism
$R_{r_\X} \cong S x_0^{r_\X}$, respectively. Then we calculate
$$
\Psi(F x_0 dG) \;=\; \phi(F) x_0 d(\phi(G)) \;=\; f x_0^{r_\X+1} \, d(g x_0^{r_\X})
\;=\; x_0^{2r_\X+1} \otimes f dg +  r_\X f g x_0^{2r_\X} dx_0
$$
Since we know already that the second summand is in the image of~$\Psi$, also the first
one is and the proof is complete.

Finally, claim (c) is an immediate consequence of~(b) and Proposition~\ref{prop:OmegaRtilde}.c,
and claim~(d) follows from Proposition~\ref{prop:TorsionMod}.a and the fact that~$\Phi$
is identified with the canonical map from $\OmegaR$ to its localization in~$x_0$.
\end{proof}

\bigbreak
%%%%%%%%%%%%%%%%%%%%%%%%%%%%%%%%%%%%%%%%%%%%%%%%%%%%%%%%
%
%  Section 4: The Euler Form and the Koszul Submodule
%
%%%%%%%%%%%%%%%%%%%%%%%%%%%%%%%%%%%%%%%%%%%%%%%%%%%%%%%%

\section{The Euler Form and the Koszul Submodule}\label{sec:EulerKoszul}

Continuing to use the notation introduced in the preceding sections, 
let~$K$ be a perfect field, let $P=K[X_0,\dots,X_n]$, and let $\X$ 
be a 0-dimensional subscheme of~$\mathbb{P}^n$ 
with homogeneous coordinate ring $R = P/I_\X$. For $i=0,\dots,n$,
we denote the image of~$X_i$ in~$R$ by~$x_i$, and we assume that~$x_0$
is a non-zerodivisor in~$R$.

Recall that the $K$-linear map $\delta:\; R \longrightarrow R$ 
given by $\delta(f) = i \cdot f$ for $i\ge 0$ and $f\in R_i$ is a derivation
of~$R$ (cf.~\cite{Kun1986}, Example~1.5). It is called the {\bf Euler derivation} of~$R/K$.

\begin{definition}\label{Euler-Form}
By the universal property of~$\OmegaR$ (cf.~\cite{Kun1986}, Theorem 1.19), the Euler
derivation $\delta: R \longrightarrow R$ gives rise to an $R$-linear map
$\epsilon:\; \OmegaR \longrightarrow R$ such that $\delta = \epsilon\circ d_{R/K}$.
This map will be called the {\bf Euler form} on~$\OmegaR$.

For $i=0,\dots,n$, we obtain $\epsilon(dx_i) = \delta(x_i) = x_i$. 
Hence we have $\Im(\epsilon) = \m$, 
where $\m=\langle x_0,...,x_n\rangle$ is the homogeneous maximal ideal of $R$,
and we will view $\epsilon$ as the surjective $R$-linear map
$$
\epsilon:\; \OmegaR \;\longrightarrow\; \m \hbox{\quad \rm given by\quad}
\epsilon(dx_i) = x_i \hbox{\quad\rm for \quad} i=0,\dots,n. 
$$
\end{definition}

The following proposition collects some properties of the kernel of the Euler form.

\begin{proposition}\label{prop:EulerProps}
Let $\epsilon:\; \OmegaR \longrightarrow \m$ be the Euler form on~$\OmegaR$.
\begin{enumerate}
\item[(a)] We have $T\OmegaR  \subseteq \Ker(\epsilon)$.

\item[(b)] We have $\Ann_R(dx_0) = \{ 0\}$ and $\Ker(\epsilon) \cap Rdx_0 = \{0\}$.

\item[(c)] For $i\ge 2 r_\X+1$, we have $\HF_{\Ker(\epsilon)}(i) = \dim_K (\OmegaS)$.

\end{enumerate}
\end{proposition}

\begin{proof}
First we show~(a). Let $w\in T\OmegaR$. By Proposition~\ref{prop:TorsionMod}.a,
there exists a number $i\ge 1$ such that $x_0^i \, w = 0$.
Since the set $\{dx_0,\dots,dx_n\}$ generates $\OmegaR$, we can write
$w = f_0 dx_0 + \cdots + f_n dx_n$ with $f_0,\dots,f_n \in R$
and calculate $0 = \epsilon (x_0^i w) = \epsilon (x_0^i f_0 dx_0 + \cdots + 
x_0^i f_n dx_n) = x_0^i f_0 x_0 + \cdots + x_0^i f_n x_n$.
Now we use the fact that~$x_0$ is a non-zerodivisor for~$R$ and conclude
that $\epsilon(w) = f_0 x_0 + \cdots + f_n x_n = 0$.

To prove~(b), we let $f\in R$ such that $fdx_0 = 0$. Then $0 = \epsilon(fdx_0) = 
f\, \epsilon(dx_0) = f\, x_0$ implies $f=0$. This shows $\Ann_R(dx_0) = \{0\}$.
Next, let $f\in R$ such that $f dx_0 \in \Ker(\epsilon)$. Then $\epsilon(fdx_0) = 
f\, x_0 = 0$ implies $f=0$, and hence $fdx_0 = 0$.

Finally, claim~(c) follows from Proposition~\ref{prop:CompareKahlerMod}.c
and $\dim_K(\m_i) = \dim_K(R_i) = \deg(\X)$ for $i \ge r_\X + 1$.
\end{proof}

Given a linear form on a module, one can construct an associated Koszul complex.

\begin{definition}\label{Euler-KoszulComplex}
As above, let~$R$ be the homogeneous coordinate ring of a 0-dimensional scheme~$\X$
in~$\mathbb{P}^n$, and let $\epsilon:\, \OmegaR \longrightarrow \m$ be the Euler form.
\begin{enumerate}
\item[(a)] The complex of $R$-modules
$$
0 \;\longrightarrow\; \Lambda^{n+1} \OmegaR \;\TTo{\epsilon^{(n+1)}}\;
\cdots \; \TTo{\epsilon^{(3)}}\;  \Lambda^2 \OmegaR  \; \TTo{\epsilon^{(2)}}\;
\OmegaR \;\TTo{\epsilon}\; \m \;\To\; 0 
$$
is called the {\bf Euler-Koszul complex} of~$R$ (or of~$\X$). Here $\epsilon^{(i)}$
is given by $\epsilon^{(i)}(w_1 \wedge \cdots \wedge w_i) = \sum_{j=1}^i
(-1)^{j-1} \epsilon(w_j)\, w_1 \wedge \cdots \wedge \widehat{w_j} \wedge
\cdots \wedge w_i$ for $i\ge 2$ and $w_1,\dots,w_i \in \OmegaR$.

\item[(b)] The image of the $R$-linear map $\epsilon^{(2)}:\,
\Lambda^2 \OmegaR \longrightarrow \OmegaR$, i.e., the $R$-sub\-module
$U_{R/K} = \langle x_i dx_j - x_j dx_i \mid 0\le i<j\le n \rangle$
of~$\OmegaR$ is called the {\bf Koszul submodule} of~$\OmegaR$.
\end{enumerate}
\end{definition}

It is clear that the Koszul submodule of~$\OmegaR$ is contained in
the kernel of the Euler form, and usually both submodules are equal.
A more precise description of their relationship is given as follows.

\begin{proposition}\label{prop:KoszulEuler}
Let $\epsilon:\, \OmegaR \longrightarrow \m$ be the Euler form and
$U_{R/K} \subseteq \OmegaR$ the Koszul submodule of $\OmegaR$.
\begin{enumerate}
\item[(a)] We have $U_{R/K} \subseteq \Ker(\epsilon)$, and this
is an equality if we have $\charac(K) = 0$ or $\charac(K) \ge 2 r_\X + 1$.

\item[(b)] Assume that $\charac(K) = p > 0$. If $i\ge 1$ satisfies $p \nmid i$ 
or if $i\ge 2 r_\X+1$ then we have $(U_{R/K})_i = \Ker(\epsilon)_i$.

\item[(c)] $\Ker(\epsilon) = \{ w \in \OmegaR \mid x_0 w \in U_{R/K} \}$.
\end{enumerate}
\end{proposition}

\begin{proof}
First of all, note that $\epsilon(x_i dx_j - x_j dx_i) = x_i x_j - x_j x_i = 0$
for $0 \le i<j \le n$. Thus we have $U_{R/K} \subseteq \Ker(\epsilon)$. Moreover,
notice that both $U_{R/K}$ and $\Ker(\epsilon)$ are graded $R$-submodules of~$\OmegaR$.

Now let $i\ge 1$. For a homogeneous element $w\in (\OmegaR)_i$, it is shown 
in the proof of~\cite{KLL2019}, Proposition~2.4, that 
$(\epsilon^{(2)}\circ d + d\circ \epsilon)(w) = i\, w$. Therefore every
element $w\in \Ker(\epsilon)_i$ satisfies $\epsilon^{(2)}(dw) = i\, w$. 
Hence if~$i$ is a unit in~$K$ then we have $w = \epsilon^{(2)}(\frac{1}{i}\, dw)
\in U_{R/K}$. This proves the case $\charac(K)=0$ in~(a) 
and the first case $p \nmid i$ in~(b).

Next we consider the second case in~(a) and the second case in~(b) simultaneously. 
By what we have already seen, we may assume that $\charac(K) = p > 0$ 
and that~$p$ divides the degree
$i\ge 2 r_\X+1$ of an element $w\in \Ker(\epsilon)_i \subseteq (\OmegaR)_i$. 
We use the fact that $\mu_{x_0}:\; (\OmegaR)_{i-1} \To (\OmegaR)_i$
is surjective by Proposition~\ref{prop:HKDprops}.b and write $w = x_0 \tilde{w}$
with $\tilde{w}\in (\OmegaR)_{i-1}$. Here $0=\epsilon(w) = x_0 \epsilon(\tilde{w})$
implies $\tilde{w}\in \Ker(\epsilon)_{i-1}$. Since $p \nmid (i-1)$, we have
$\tilde{w} = \frac{1}{i-1}\, \epsilon^{(2)}(d\tilde{w}) \in U_{R/K}$, and thus
$w = x_0 \tilde{w} \in U_{R/K}$.
  
It remains to prove~(c). By~(b), the only case where the containment
$(U_{R/K})_i \subseteq \Ker(\epsilon)_i$ can be strict for some $i\ge 1$ is when
$\charac(K)=p>0$ and $p\mid i$. Then $p \nmid (i+1)$ implies $x_0w \in (U_{R/K})_{i+1}$.
This proves the inclusion $\subseteq$ of the claim. 
Con\-versely, if $x_0w \in U_{R/K}\subseteq \Ker(\epsilon)$, then $0 = \epsilon(x_0w)
= x_0 \epsilon(w)$ shows $\epsilon(w) =0$.
\end{proof}

Let us also try to describe those elements in $\Ker(\epsilon)$
which are not in the Koszul submodule of $\OmegaR$ as explicitly as possible.
The following construction will be useful.

\begin{remark}
Let $P = K[X_0,\dots,X_n]$ and $\M = \langle X_0,\dots,X_n\rangle$ its 
homogeneous maximal ideal. Then every polynomial $F\in \M$ has a unique decomposition
of the form $F= F_0 X_0 + \cdots + F_n X_n$ where $F_i \in K[X_i,\dots,X_n]$. In other words,
we have $F_n \in K[X_n]$, $F_{n-1} \in K[X_{n-1},X_n]$, etc.
We call it the {\bf triangular decomposition} of~$F$. 

If we denote the tuple $(F_0,\dots,F_n)$ by~$\theta(F)$,
then the map $\theta:\, \M \To P^{n+1}$ is $K$-linear and yields a section
of the epimorphism $\xi:\, P^{n+1} \To \M$ given by $\xi(e_i)=X_i$ for $i=0,\dots,n$.
In other words, the map $\theta$ splits the short exact sequence
$$
0 \;\To\; \Ker(\xi) \;\To\; P^{n+1} \; \TTo{\xi}\; \M \;\To\; 0
$$
Recall that $\Ker(\xi)= \langle X_i e_j - X_j e_i \mid 0 \le i<j\le n\rangle$,
as $(X_0,\dots,X_n)$ is a $P$-regular sequence.
\end{remark}

With the help of the triangular decomposition, we obtain the following 
description of the kernel of the Euler form.

\begin{proposition}{\bf (The Kernel of the Euler Form)}\label{prop:KerEuler}\\
Let $R=P/I_\X$ be the homogeneous coordinate ring of a 0-dimensional subscheme~$\X$
of~$\mathbb{P}^n$, and let $\epsilon:\, \OmegaR \To \m$ be the Euler form.
\begin{enumerate}
\item[(a)] Let $\gamma:\, P^{n+1} \To \OmegaR$ be the epimorphism
defined by $\gamma(e_i) = dx_i$ for $i=0,\dots,n$, and let $\eta:\, \OmegaR \To
\OmegaR / U_{R/K}$ be the canonical epimorphism. Then the composition
$\delta = \eta\circ \gamma \circ \theta:\, \M \To \OmegaR / U_{R/K}$ is a $P$-linear map.

\item[(b)] We have $\delta(I_\X) = \Ker(\epsilon)/U_{R/K}$. 

\item[(c)] Given a homogeneous system of generators $\{G_1,\dots,G_r\}$
of~$I_\X$, we have $\Ker(\epsilon) = U_{R/K} + \langle (\gamma\circ\theta)(G_1),\dots,
(\gamma\circ\theta)(G_r)\rangle$.

\item[(d)] Let $m_\X = \max \{ i\ge 0 \mid (I_\X/ \M I_\X)_i \ne 0 \}$ be the 
maximal degree of a homogeneous minimal generator of~$I_\X$. Then we have
$(U_{R/K})_i = \Ker(\epsilon)_i$ for all $i > m_\X$.
\end{enumerate}
\end{proposition}

\begin{proof}
First we prove~(a). As a composition of $K$-linear maps, the map $\delta$
is clearly $K$-linear. Now let $F\in\M$, and let $F=F_0 X_0 + \cdots + F_n X_n$
be the triangular decomposition of~$F$, where $F_i \in K[X_i,\dots,X_n]$.

Given a term $T\in K[X_0,\dots,X_n]$, we let~$k$ be the least index of an indeterminate
dividing~$T$, and we write $T=X_k \cdot T'$ with $T' \in \mathbb{T}^{n+1}$. Then we
have
$$
T\, F \;=\;  \tsum_{i=0}^{k-1} T\, F_i\, X_i \;+\; (\tsum_{i=k}^n F_i X_i T')\, X_k
$$
and this is clearly the triangular decomposition of $T\, F$. Hence we get
$$
\delta(T\, F) \;=\; \tsum_{i=0}^{k-1} t\, f_i\, dx_i + \tsum_{i=k}^n f_i\, x_i\, t'\, dx_k
+ U_{R/K}
$$
where $f_i$ denotes the residue class of~$F_i$ in~$R$ for $i=0,\dots,n$ and~$t$
(resp.~$t'$) the residue class of~$T$ (resp.~$T'$). On the other hand, we calculate
\begin{align*}
t \cdot \delta(F) & \;=\; \tsum_{i=0}^{k-1} t\, f_i\, dx_i + \tsum_{i=k}^n x_k t'\, f_i\, dx_i 
+ U_{R/K}\\ 
& \;=\; \tsum_{i=0}^{k-1} t\, f_i\, dx_i + \tsum_{i=k}^n t' f_i x_i dx_k + \tsum_{i=k}^n
t' f_i (x_k dx_i - x_i dx_k) + U_{R/K} \\
& \;=\; \delta(T\,F)
\end{align*}
Consequently, the map~$\delta$ is $P$-linear.

To show~(b), we let $F\in I_\X$ and write the triangular decomposition of~$F$
as above as $F = F_0 X_0 + \cdots + F_n X_n$. Then we get $\delta(F) = 
f_0 dx_0 + \cdots + f_n dx_n + U_{R/K}$ and for every element $u\in U_{R/K}$
we obtain $\epsilon(f_0 dx_0 + \cdots + f_n dx_n + u) = f_0 x_0 + \cdots + f_n x_n 
= \allowbreak F + I_\X = 0$. 

For the reverse inclusion, let 
$w \in \Ker(\epsilon)$, and write $w= g_0 dx_0 + \cdots + g_n dx_n$. Represent 
$g_i\in R$ by polynomials $G_i \in P$, and let 
$F = G_0 X_0 + \cdots + G_n X_n$.
Since $\epsilon(w) = g_0 x_0 + \cdots + g_n x_n = 0$, we have $F \in I_\X$.
Now let $F= F_0 X_0 + \cdots + F_n X_n$ be the triangular decomposition of~$F$.
Then $(F_0-G_0,\dots,F_n-G_n)$ is a syzygy of $(X_0,\dots,X_n)$, and thus contained
in $V := \langle X_i e_j -X_j e_i \mid 0\le i<j\le n\rangle$. This implies
$$
\theta(F) = (F_0,\dots,F_n) = (G_0,\dots,G_n) + v \hbox{\qquad \rm with\;}v\in V
$$
and therefore $\delta(F) = g_0 dx_0 + \cdots + g_n dx_n + u + U_{R/K} = w + U_{R/K}$,
where $u= \gamma(v) \in U_{R/K}$.

Since claim~(c) follows immediately from~(a) and~(b), it remains to prove~(d).
If we have $(U_{R/K})_{m_\X} = \Ker(\epsilon)_{m_\X}$ then~(c) shows that 
the claim holds. So, the only remaining case is $\charac(K) = p>0$ 
and $p \mid m_\X$. In this case we have $(U_{R/K})_{m_\X +1} = \Ker(\epsilon)_{m_\X +1}$ 
by Proposition~\ref{prop:KoszulEuler}.a,b. Hence we see that 
$f\cdot (\gamma\circ\theta)(G_i) \in U_{R/K}$ for all $i\in \{1,\dots,r\}$
and all homogeneous elements $f\in R$ such that $\deg(f) + \deg(G_j) = m_\X+1$.
The claim follows from~(c) and this observation.
\end{proof}

The next example shows that the torsion submodule of~$\OmegaR$ is, in general, not
contained in its Koszul submodule.

\begin{example}
Let $K=\mathbb{F}_3$, let $P=K[x_0,x_1,x_3]$, and let~$\X$ be the 0-dimensional 
complete intersection scheme in $\mathbb{P}^2$ defined by 
the homogeneous ideal $I_\X = \langle G_1,G_2\rangle$, where
$G_1 = X_1^2 + X_2^2$ and $G_2 = X_0 X_1^2 + X_1^3 + X_2^3$.

It is straightforward to calculate $\HF_\X:\ 1\ 3\ 5\ 6\ 6\cdots$ which shows $r_\X=3$,
as well as $\HF_{\OmegaR}:\ 0\ 3\ 8\ 11\ 10\ 10\cdots$ which implies $\ri(\OmegaR) = 4$.
Moreover, we compute
\begin{align*}
&\HF_{T\OmegaR}:\ 0\ 0\ 0\ 1\ 0\ 0\cdots, \hbox{\ \rm and thus\ } \ri(T\OmegaR) = 4,\\
&\HF_{U_{R/K}}: \,\,\  0\ 0\ 3\ 4\ 4\ 4\cdots,\hbox{\ \rm and thus\ } \ri(U_{R/K}) = 3,\\
&\HF_{\Ker(\epsilon)}: \,\ 0\ 0\ 3\ 5\ 4\ 4\cdots, \hbox{\ \rm and thus\ } \ri(\Ker(\gamma)) = 4.
\end{align*}
Consequently, the Koszul submodule $U_{R/K}$ is a proper submodule of~$\Ker(\epsilon)$.

The triangular decomposition of $G_2\in I_\X$ is 
$G_2 = (X_{0}X_{1}+ X_{1}^{2})\, X_1 + (X_2^2)\, X_2$.
Therefore we may check that
$$
(\gamma\circ\theta)(G_2) \;=\; (x_0 x_1 + x_1^2)\, dx_1 + (x_2^2)\, dx_2 \in \Ker(\epsilon)
\setminus U_{R/K}
$$
Furthermore, it turns out that $T\OmegaR = K\cdot (\gamma\circ\theta)(G_2)$, and hence
the torsion submodule of~$\OmegaR$ is in general not contained in its Koszul submodule.
Notice that for reduced schemes we have $T\OmegaR = \Ker(\epsilon)$ 
(see~\cite{DK1999}, Proposition 2.1). However, even in this case 
$U_{R/K} \subsetneq \Ker(\epsilon)$ may occur, for instance, when $\X$ is a set of 
five $\mathbb{F}_3$-rational points on a non-singular conic in $\mathbb{P}^2$.
\end{example}

\bigbreak
%%%%%%%%%%%%%%%%%%%%%%%%%%%%%%%%%%%%%%%%%%%%%%%%%%%%
%
%  Section 5: Kaehler Differential $m$-Forms
%
%%%%%%%%%%%%%%%%%%%%%%%%%%%%%%%%%%%%%%%%%%%%%%%%%%%%

\section{K\"ahler Differential $m$-Forms}\label{sec4:KDmForms}

Continuing to use the above notation and assumptions,
we want to introduce and study the higher exterior powers of~$\OmegaR$
next.

\begin{definition}
Let $m\ge 0$, and let $R$ be the homogeneous coordinate ring of a zero-dimensional
scheme~$\X$ in~$\mathbb{P}^n$.
\begin{enumerate}
\item[(a)] The $m$-th exterior power $\OmegaRm := \Lambda^m_R \OmegaR$
is called the {\bf module of K\"ahler differential $m$-forms} of~$R$.

\item[(b)] The graded $R$-algebra $\Omega^{\bullet}_{R/K} :=
\bigoplus_{m\ge 0} \OmegaRm$ is called the {\bf K\"ahler differential algebra} of $R/K$.
\end{enumerate}
\end{definition}

Notice that $\Omega^0_{R/K} = R$ and that we can define $P$-modules
$\Omega^m_{P/K}$ analogously. For every $m\ge 0$, the module of K\"ahler differential
$m$-forms is a finitely generated graded $R$-module. To calculate a presentation
of this $R$-modules, we can use the following results.

\begin{proposition}\label{prop:PresentationOmegaRm}
Let $R$ be the homogeneous coordinate ring of a zero-dimensional
scheme~$\X$ in~$\mathbb{P}^n$.
\begin{enumerate}
\item[(a)] For $m\ge n+2$, we have $\OmegaRm = \langle 0 \rangle$.

\item[(b)] For $m=0,\dots,n+1$, the $P$-module $\Omega^m_{P/K}$ is a free
$P$-module of rank $\binom{n+1}{m}$ with basis
$\{ dX_{i_1} \wedge \cdots \wedge
dX_{i_m} \mid 0\le i_1 < \cdots < i_m \le n \}$.

\item[(c)] For $m\ge 2$, the module of K\"ahler differential $m$-forms
satisfies
$$
\OmegaRm \;\cong\;   \Omega^m_{P/K} \;/\; \left(  I_\X \Omega^m_{P/K} +
dI_\X \wedge \Omega^{m-1}_{P/K}  \right)
$$
Here $dI_\X$ is the $P$-submodule of~$\Omega^1_{P/K}$ generated by
$\{d_{P/K}f \mid f\in I_\X\}$.

\item[(d)] Let $I_\X = \langle f_1,\dots,f_k \rangle$ for some $f_1,\dots,f_k\in P$.
Then we have
$$
d I_\X = \langle df_1,\dots,df_k \rangle + \langle f_i dX_j
\mid i\in \{1,\dots,k\},\; j\in \{0,\dots,n\} \rangle
$$
\end{enumerate}
\end{proposition}

\begin{proof}
Claim (a) is clearly true, claim (b) follows from \cite{Kun1986}, Example 2.5.d,
and claim~(c) follows from ibid., Proposition~4.12.

To show~(d), we note that $I_{\X}\Omega^1_{P/K}$ is contained in $dI_\X$
because, for $f\in I_\X$ and $g\in P$, we have $f\, dg = d(fg) - g\, df
\in dI_\X$. This implies the inclusion $\supseteq$. Conversely, let
$g\in I_\X$ and write $g=h_1 f_1 + \cdots + h_k f_k$ with $h_1,\dots,h_k \in P$.
Then $dg = h_1 df_1 + \cdots + h_k df_k + f_1 dh_1 + \cdots + f_k dh_k$
is contained in the right-hand side, and the desired equality is proved.
\end{proof}

Since the modules of K\"ahler differential $m$-forms are
finitely generated graded $R$-modules, we can define their
Hilbert function, Hilbert polynomial and regularity index as usual.

\begin{definition}
Let $R$ be the homogeneous coordinate ring of a zero-dimensional
scheme~$\X$ in~$\mathbb{P}^n$, and let $m\ge 1$.
\begin{enumerate}
\item[(a)] The map $\HF_\OmegaRm:\; \ZZ \To \ZZ$ given by
$\HF_\OmegaRm (i) = \dim_K (\OmegaRm)_i$ for every $i\in \ZZ$ is called
the {\bf Hilbert function} of~$\OmegaRm$.

\item[(b)] The integer $\HP(\OmegaRm) = \HF_{\OmegaRm} (i)$ for $i\gg 0$
is called the {\bf Hilbert polynomial} of~$\OmegaRm$.

\item[(c)] The number $\ri(\OmegaRm) := \min \{ i\in\ZZ \mid
\HF_{\OmegaRm}(j)=\HP(\OmegaRm) \hbox{\ \rm for\ }j\ge i\}$
is called the {\bf regularity index} of~$\OmegaRm$.
\end{enumerate}
\end{definition}

Our goal in this section is to determine the Hilbert polynomial of the
module~$\OmegaRm$ and to find a sharp bound for its regularity index.
As in the preceding section, our method to reach these goals is to compare
it to the module of K\"ahler differential $m$-forms of the
truncated integral closure~$\widetilde{R}$ of~$R$ in~$Q^h(R)$.
The modules $\OmegaRtm$ can be computed and compared to $\OmegaRm$
as follows.

\begin{proposition}\label{prop:HigherKDM}
Let $R$ be the homogeneous coordinate ring of a zero-dimensional
scheme~$\X$ in~$\mathbb{P}^n$, let $\widetilde{R}$ be the truncated integral
closure of~$R$ in~$Q^h(R)$, let $S=R/\langle x_0-1\rangle$ be the
affine coordinate ring of~$\X$ in~$\mathbb{A}^n \cong D_+(x_0)$,
and let $m\ge 1$.
\begin{enumerate}
\item[(a)] The isomorphism $\widetilde{R} \cong S[x_0]$ induces
an isomorphism
$$
\Psi^{(m)}:\; \OmegaRtm  \;\cong\; K[x_0] \otimes_K \Omega^m_{S/K} \;\oplus\;
K[x_0]\, dx_0 \wedge  \Omega^{m-1}_{S/K}
$$

\item[(b)] 
Let $\Phi:\; \OmegaR \To \OmegaRtilde$ be the canonical $R$-linear map
induced by $R \subseteq \widetilde{R}$. Then the $K$-linear map
$(\Lambda^m_R \Phi)_i:\;  (\OmegaRm)_i \To (\OmegaRtm)_i$ is an isomorphism for
every $i\ge 2 r_\X+m$.

\item[(c)] We have $\HP(\OmegaRm) = \HP(\OmegaRtm) = \dim_K(\Omega^m_{S/K}) +
\dim_K(\Omega^{m-1}_{S/K})$.

\item[(d)] We have $\ri(\OmegaRm) \le 2 r_\X + m$.

\item[(e)] If $\charac(K) = 0$ or if $\charac(K)=p>0$ and $p \nmid (2r_\X+n+1)$, 
then we even have $\ri(\Omega^{n+1}_{R/K})\le 2 r_\X+n$.
\end{enumerate}
\end{proposition}

\begin{proof}
Claim (a) follows from
\begin{align*}
	\OmegaRtm &= \Lambda^m_{\widetilde{R}}
	\;\OmegaRtilde \cong \Lambda^m_{S[x_0]}\;
	\left( (S\otimes_K \Omega^1_{K[x_0]/K}) \oplus ( K[x_0]
	\otimes_K \Omega^1_{S/K}) \right)\\
	&\cong {\textstyle\bigoplus\limits}_{i=0}^{m}
	\Lambda^i_{S[x_0]} (S \otimes_K\Omega^1_{K[x_0]/K} ) \wedge
	\Lambda^{m-i}_{S[x_0]} (K[x_0]  \otimes_K \OmegaS)\\
	&\cong K[x_0] \otimes_K \Omega^m_{S/K} \oplus (S \otimes_K \Omega^1_{K[x_0]/K})
	\wedge ( K[x_0] \otimes_K \Omega^{m-1}_{S/K} )
\end{align*}
To prove (b), we first show that $\widetilde{\Phi} := (\Lambda^m_R \Phi)_i$
is surjective for $i\ge 2r_\X+m$. Notice that
$K[x_0] \otimes_K \Omega^m_{S/K} \cong S[x_0]\otimes_S \Omega^m_{S/K}$.
Using (a), it suffices to show that
$$
x_0^{i-1}dx_0\wedge f dx_{i_1}\wedge \cdots \wedge dx_{i_{m-1}},\;
x_0^{i}\otimes g dx_{j_1}\wedge \cdots \wedge dx_{j_m} \in \Im(\widetilde{\Phi})
$$
where $f,g\in S$ and
$\{i_1,...,i_{m-1}\}, \{j_1,...,j_{m}\}$ are subsets of $\{1,...,n\}$
of lengths $m-1$ and $m$.
Let $F,G\in R_{r_\X}$ be preimages of $f x_0^{r_\X}$ and
$g x_0^{r_\X}$ under the isomorphism $R_{r_\X} \cong S x_0^{r_\X}$ 
in Proposition~\ref{prop:RandRtilde}.b, respectively.
Then we have
$$
\begin{aligned}
\widetilde{\Phi} &(Fx_0^{i-r_\X-m}dx_0\wedge dx_{i_1}\wedge \cdots \wedge dx_{i_{m-1}})\\
&=\Phi(Fx_0^{i-r_\X-m}dx_0)\wedge \Phi(dx_{i_1})
\wedge\cdots\wedge \Phi(dx_{i_{m-1}})\\
&=  F\deh x_0^{i-m}dx_0 \wedge (x_{i_1}dx_0 + x_0\otimes dx_{i_1})
\wedge \cdots\wedge (x_{i_{m-1}}dx_0 + x_0\otimes dx_{i_{m-1}})\\
&= x_0^{i-1}dx_0\wedge f dx_{i_1}\wedge \cdots \wedge dx_{i_{m-1}}
\end{aligned}
$$
Also
$$
\begin{aligned}
\widetilde{\Phi}&(Gx_0^{i-r_\X-m} dx_{j_1}\wedge dx_{j_2}\wedge \cdots \wedge dx_{j_{m}})\\
	&=\Phi(Gx_0^{i-r_\X-m} dx_{j_1})\wedge \Phi(dx_{j_2})
	\wedge\cdots\wedge \Phi(dx_{j_{m}})\\
	&=  (G\deh x_{j_1} x_0^{i-m}dx_0 + x_0^{i-m+1}\otimes G\deh dx_{j_1}) \wedge \\
	&\quad \wedge (x_{j_2}dx_0 + x_0\otimes dx_{j_2})
    \wedge \cdots\wedge (x_{j_m}dx_0 + x_0\otimes dx_{j_m})\\
	&= x_0^{i-1}dx_0\wedge g\textstyle{\sum\limits_{k=1}^m}x_{j_k} dx_{j_1}
	\wedge \cdots \wedge \widehat{dx_{j_k}}\wedge \cdots \wedge dx_{j_{m}} +\\
	&\quad + x_0^{i}\otimes gdx_{j_1}\wedge dx_{j_2}\wedge \cdots \wedge dx_{j_{m}}
\end{aligned}
$$
where $\widehat{dx_{j_k}}$ means that $dx_{j_k}$ is omitted.
In the last equality, the first summand is in the image of~$\widetilde{\Phi}$,
and so is the second one. Therefore $\widetilde{\Phi}$ is surjective.
Next, we show that $\dim_K(\OmegaRm)_i=\dim_K(\OmegaRtm)_i$
for $i\ge 2r_\X+m$, and hence $\widetilde{\Phi}$ is bijective.
Set $W:=\Im(\Phi)\subseteq \OmegaRtilde$.
Then $\Phi: \OmegaR \rightarrow W$ is an epimorphism
of graded $R$-modules. By \cite{SS1988}, Satz~85.8, we have
the exact sequence of $K$-vector spaces
$$
0  \;\longrightarrow\;  (\Ker(\Phi)\wedge_R \Omega^{m-1}_{R/K})_i
\;\longrightarrow\; (\Omega^{m}_{R/K})_i  \;\stackrel{\widetilde{\Phi}}{\longrightarrow}\;
(\Lambda_R^{m} (W))_i   \;\longrightarrow\;   0.
$$
It follows from Proposition~\ref{prop:CompareKahlerMod}.d 
that $T\OmegaR = \Ker(\Phi)$. By Proposition~\ref{prop:TorsionMod},  
$\Ker(\Phi)_k = \{0\}$ for all $k\ge 2r_\X+1$, and so
$(\Ker(\Phi)\wedge_R \Omega^{m-1}_{R/K})_i = \{0\}$.
It follows that  $\dim_K(\Omega^{m}_{R/K})_i=\dim_K(\Lambda_R^{m} (W))_i$.
On the other hand, it is clearly true that
$(\Lambda_R^{m} (W))_i\subseteq (\OmegaRtm)_i$.
Hence the surjectivity of $\widetilde{\Phi}$ implies that
$(\Lambda_R^m (W))_i = (\OmegaRtm)_i$,
and subsequently we get $\dim_K(\OmegaRm)_i=\dim_K(\OmegaRtm)_i$,
as desired.

Claims (c) and~(d) follow from~(a) and~(b). It remains to prove~(e).
Consider the Euler-Koszul complex of~$R$
$$
0 \;\longrightarrow\; \Lambda^{n+1} \OmegaR \;\TTo{\epsilon^{(n+1)}}\;
\cdots \; \TTo{\epsilon^{(3)}}\;  \Lambda^2 \OmegaR  \; \TTo{\epsilon^{(2)}}\;
\OmegaR \;\TTo{\epsilon}\; \m \;\To\; 0 
$$
This sequence is an exact sequence of graded $R$-modules in characteristic zero
by~\cite{KLL2019}, Proposition~2.4. If $\charac(K)=p>0$ and $p \nmid
(2 r_\X +n+1)$ then the proof of \cite{KLL2019}, Proposition~2.4, shows that the
map $\epsilon^{(n+1)}:\; (\Omega^{n+1}_{R/K})_d \longrightarrow (\Omega^n_{R/K})_d$
is also injective in degree $d= 2 r_\X + n+1$. In both cases we see that 
the claim follows from~(d).
\end{proof}

The above upper bound for the regularity index of $\OmegaRm$
is sharp, as the following example shows.

\begin{example}
Let $\X=\{p_1,...,p_{n+1}\}$ be a set of $n+1$ distinct $K$-rational
points in~$\mathbb{P}^n$ which are in general position. We claim that
the regularity index of $\OmegaRm$
satisfies $\ri(\OmegaRm)=m+2$ for $1\le m\le n$
and $\ri(\Omega^{n+1}_{R/K})=n+2$.

First of all, it is clearly true that $\HF_\X:\ 1\ n+1\ n+1\cdots$
and $r_\X =1$. After a linear change of coordinates, we may assume that
$p_1=(1:0:\dots:0)$, $\dots$, $p_{n+1}=(0:\dots:0:1)$.
It follows that $I_\X = \langle X_i X_j \mid 0\le i<j\le n\rangle$.

Now let~$\Y$ be the 0-dimensional subscheme of~$\mathbb{P}^n$ defined by the saturated 
homogeneous ideal $I_\Y = \bigcap_{j=1}^{n+1} I_{p_j}^2$. Then we have
$\HF_\Y: 1\;\binom{n+1}{n}\; \binom{n+2}{n}\; (n+1)^2\; (n+1)^2\cdots$ and $r_\Y = 3$.
From \cite{Kun1986}, Proposition~4.17, we obtain an exact sequence
of graded $R$-modules
$$
0 \;\longrightarrow\; I_\X/I_\Y  \;\longrightarrow\; 
\Omega^1_{P/K}/I_\X \Omega^1_{P/K} \;\longrightarrow\; \OmegaR 
\;\longrightarrow\; 0
$$
and we deduce that $\ri(\OmegaR)=3$. Moreover,
$\Omega^{n+1}_{R/K} \cong (P/\langle X_0,...,X_n\rangle)(-n-1)$,
implies $\HF_{\Omega^{n+1}_{R/K}}:\ 0 \dots 0\ 1\ 0\ 0\dots$
and $\ri(\Omega^{n+1}_{R/K})=n+2$.

Altogether, by Proposition~\ref{prop:SmoothnessCriterion}, it 
suffices to show $\HF_{\Omega^{m}_{R/K}}(m+1)\ne 0$ for all $m\in \{2,\dots,n\}$.
Observe that the element $w = dx_0\wedge dx_1\wedge\cdots\wedge dx_n$
in $(\Omega^{n+1}_{R/K})_{n+1}$ is non-zero.
Since $w = d(x_0\, dx_1\wedge\cdots\wedge dx_n) \ne 0$, also the element 
$\tilde{w} = x_0\, dx_1\wedge\cdots\wedge dx_n$ in $\Omega^{n}_{R/K}$ is non-zero.
This implies $\HF_{\Omega^n_{R/K}}(n+1)\ne 0$.  

As $\tilde{w} = (x_0\, dx_1 \wedge\cdots\wedge dx_m) \wedge (dx_{m+1} 
\wedge \cdots \wedge dx_n)$, we obtain non-zero elements
$x_0\, dx_1\wedge \cdots \wedge dx_m$ in $\OmegaRm$ for all $m\in \{2,\dots,n \}$,
and the claim is proved.
\end{example}

\begin{remark}
Let $\X\subseteq \mathbb{P}^n$ be a set of $t$ distinct $K$-rational points.
If $t\le n$ then $\Omega^{n+1}_{R/K}=0$.
When $t>n$ and $r:=\ri(\Omega^{n+1}_{R/K})>0$,
we get a lower bound for the regularity index of $\OmegaRm$ given by
$\ri(\OmegaRm)\ge r+m-n$ for $2\le m\le n$. In particular,
if $r=2r_\X+n$ then $\ri(\OmegaRm)= 2r_\X+m$ for $2\le m\le n$.
\end{remark}

\bigbreak
%%%%%%%%%%%%%%%%%%%%%%%%%%%%%%%%%%%%%%%%%%%%%%%%%%%%%%%%%%%%
%
%  Section 6: K\"ahler Differentials of Curvilinear Schemes
%
%%%%%%%%%%%%%%%%%%%%%%%%%%%%%%%%%%%%%%%%%%%%%%%%%%%%%%%%%%%%

\section{K\"ahler Differentials of Curvilinear Schemes}\label{sec5:CurvilinearZD}

As before, we let $\X$ be a 0-dimensional scheme in~$\mathbb{P}^n$
over a perfect field~$K$, we let $R=P/I_\X$ be the homogeneous coordinate
ring of~$\X$, we assume that $x_0\in R_1$ is a non-zero divisor of~$R$, and we
let $S=R/\langle x_0-1\rangle$ be the affine coordinate ring of~$\X$
in~$\mathbb{A}^n \cong D_+(x_0)$.
In this section we want to connect local properties of~$\X$ to the structure
of the modules of K\"ahler differential $m$-forms of~$R/K$ for various $m\ge 1$.

Our first proposition characterizes smooth schemes.

\begin{proposition}\label{prop:SmoothnessCriterion}
In the above setting, the following conditions are equivalent.
\begin{enumerate}
\item[(a)] The scheme~$\X$ is smooth.

\item[(b)] We have $\OmegaS = \langle 0 \rangle$.

\item[(c)] We have $\HP(\OmegaR)=\deg(\X)$ and $\HP(\OmegaRm)=0$ for $m\ge 2$.
\end{enumerate}
\end{proposition}

\begin{proof}
The equivalence of (a) and~(b) follows from \cite{Kun1986}, Corollary 6.9,
and the equivalence of~(b) and~(c) follows from Proposition~\ref{prop:CompareKahlerMod}.
\end{proof}

Recall that a 0-dimensional scheme~$\X$ in~$\mathbb{A}^n$ has
a finite support and that $S$ is a 0-dimensional affine $K$-algebra.
Let $\q_1,\dots,\q_t$ be the primary components
of $\langle 0\rangle$ in~$S$.
Then the Chinese Remainder Theorem yields an isomorphism
$S \cong S/\q_1 \times \cdots \times S/\q_t$ where the rings
$\mathcal{O}_i = S/\q_i$ are 0-dimensional local rings.

\begin{definition}\label{def:curvilinear}
In the above setting, let $(T,\n)$ be a 0-dimensional local ring
with residue class field $L=T/\n$.
\begin{enumerate}
\item[(a)] The ring~$T$ is called {\bf (weakly) curvilinear} if
its maximal ideal~$\n$ is generated by one element.

\item[(b)] The scheme~$\X$ is called {\bf (weakly) curvilinear} if
its local rings $\mathcal{O}_1,\dots, \mathcal{O}_t$ are all (weakly) 
curvilinear.
\end{enumerate}
\end{definition}

The condition to be weakly curvilinear can be characterized as follows.

\begin{proposition}\label{prop:CharCurvilinear}

Let $(T,\n)$ be a 0-dimensional local Noetherian ring
with residue class field $L=T/\n$, and let 
$\gr_{\n}(T)$  be the {\bf associated graded ring} of~$T$.
Then the following conditions are equivalent.
\begin{enumerate}
\item[(a)] The ring~$T$ is weakly curvilinear.

\item[(b)] The associated graded ring of~$T$ is of the form
$\gr_{\n}(T) \cong L[z]/\langle z^m\rangle$
for some $m\ge 1$.

\end{enumerate}
\end{proposition}

\begin{proof}
First we prove that (a) implies~(b).
Notice that the ring~$T$ is Artinian, and let $\n=\langle a\rangle$
for some $a\in T$. Then every proper ideal of~$T$ is a power of~$\n$, 
and therefore principal. In particular, there is a smallest number $m\ge 1$
such that $\langle 0\rangle = \n^m
=\langle a^m\rangle$. Thus the $L$-algebra epimorphism
$L[z]\rightarrow \gr_{\n}(T)$ defined by
$z\mapsto a + \n^2$ yields an isomorphism
$\gr_{\n}(T) \cong L[z]/\langle z^m \rangle$, as desired.

Conversely, if $\gr_{\n}(T) \cong L[z]/\langle z^m\rangle$, then
$\dim_L(\n/\n^2)=1$, and hence $\n$
is principal by Nakayama's Lemma. So, the ring~$T$ is weakly curvilinear.
\end{proof}

Let us also compare the above definition of a weakly curvilinear
0-dimensional scheme to the one given in~\cite{KR2016}.

\begin{remark}
In~\cite{KR2016}, Definition~4.3.1, a 0-dimensional local affine
$K$-algebra~$T$ is called {\bf weakly curvilinear} if it is of the form
$T \cong K[z]/\langle p(z)^m \rangle$, where $p(z)$ is an
irreducible polynomial in~$K[z]$ and $m\ge 1$, and it is called {\bf curvilinear}
if $T \cong K[z]/\langle z^m\rangle$ for some $m\ge 1$.

For a ring~$T$ which is weakly curvilinear in this sense, its maximal 
ideal $\langle p(\bar{z}) \rangle$ is a principal ideal. Therefore 
Definition~\ref{def:curvilinear} generalizes the definition in~\cite{KR2016}.
Since the definition of a weakly curvilinear 0-dimensional scheme uses
the local rings of its affine coordinate ring, also this definition
generalizes the one in~\cite{KR2016}. 

Moreover, note that in the case of an algebraically closed field~$K$,
the maximal ideal of a local ring $\mathcal{O}_i$ of~$\X$ is unigenerated
if and only if it is the residue class ideal of an ideal of the form
$\langle z-a\rangle$ in $K[z]$, where $a\in K$. Thus the ring $\mathcal{O}_i$
is isomorphic to a ring of the form $K[z]/\langle z^m\rangle$ with $m\ge 1$,
in agreement with the classical definition of a curvilinear 0-dimensional scheme.
Such schemes appear naturally in algebraic geometry, as for instance the 
well-known paper~\cite{EH1992} demonstrates.
\end{remark}

The next example shows that Definition~\ref{def:curvilinear} is a proper
generalization of Definition~4.3.1 in~\cite{KR2016}.

\begin{example}
Consider the 0-dimensional scheme $\X$ in $\mathbb{P}^2$ 
over the field~$\mathbb{Q}$ defined by the ideal 
$I_\X =\langle X_1^2+X_0^2,(X_2^2-2X_0)^2\rangle
\subseteq \mathbb{Q}[X_0,X_1,X_2]$. Its affine coordinate ring 
$S=\mathbb{Q}[X_1,X_2]/\langle X_1^2+1,(X_2^2-2)^2\rangle$ is
a 0-dimensional local ring with maximal ideal
$\n = \langle x_2^2-2 \rangle$ and residue class field 
$L = S/\n \cong \mathbb{Q}(i,\sqrt{2})$.

Since the ideal~$\n$ is principal, the ring~$S$  
is weakly curvilinear in the sense of Definition~\ref{def:curvilinear}.
However, the ring~$S$ cannot be presented as $S\cong \mathbb{Q}[z]/\langle p(z)^m\rangle$ 
with an irreducible polynomial $p(z)\in \mathbb{Q}[z]$ and $m\ge 1$.
Furthermore, note that $S \ncong \gr_{\n}(S)$.
\end{example}

The main topic of this section is to use the Hilbert polynomials of the 
modules of K\"ahler differential $m$-forms to characterize weakly curvilinear
0-dimensional schemes.

\begin{proposition}\label{prop:CharCurvilinearKD}
For a 0-dimensional scheme~$\X$ in~$\mathbb{P}^n$ as above,
the following conditions are equivalent.
\begin{enumerate}
\item[(a)] The scheme~$\X$ is weakly curvilinear, but not smooth.

\item[(b)] We have $\OmegaS \ne \langle 0 \rangle$ and $\Omega^m_{S/K} =
\langle 0 \rangle$ for all $m\ge 2$.

\item[(c)] We have $\HP(\OmegaR) > \deg(\X)$ and
$\HP(\Omega^2_{R/K}) = \HP(\OmegaR) - \deg(\X)$
and $\HP(\OmegaRm) = 0$ for $m > 2$.

\end{enumerate}
\end{proposition}

\begin{proof}
As above, we let $S= R/\langle x_0-1\rangle$ be the affine coordinate
ring of~$R$ and $S \cong \mathcal{O}_1 \times \cdots \times \mathcal{O}_t$
its decomposition into local rings.

First we show that~(a) implies~(b). Since~$\X$ is not smooth,
we have $\OmegaS \ne \langle 0\rangle$ by
Proposition~\ref{prop:SmoothnessCriterion}. Using the isomorphism
$\Omega^2_{S/K} \cong \Omega^2_{\mathcal{O}_1/K} \times \cdots
\times \Omega^2_{\mathcal{O}_t/K}$
(see \cite{Kun1986}, Proposition~4.7), we see that it suffices to
show $\Omega^2_{\mathcal{O}_i/K}= \langle 0 \rangle$ for $i=1,\dots,t$.
As~$\X$ is weakly curvilinear, the maximal ideal of~$\mathcal{O}_i$
is unigenerated. Then the differential of this element generates
the $\mathcal{O}_i$-module $\Omega^1_{\mathcal{O}_i/K}$. It follows
that the second and all higher exterior powers of this module are zero.

Conversely, let us show that~(b) implies~(a).
Since $\OmegaS \ne \langle 0 \rangle$, the scheme~$\X$ is not smooth.
For $i=1,\dots,t$, the hypothesis and the above isomorphism
yield $\Omega^2_{\mathcal{O}_i/K} = 0$.

Hence we may assume that~$S$ is local. Let
$\n$ be the maximal ideal of~$S$. Then we have
$\n / \n^2 \cong \OmegaS / \n\,\OmegaS$ by~\cite{Kun1986}, Corollary~6.5,
and this $S/\n$-vector space satisfies
$$
\Lambda^2_{S/\n} \, (\OmegaS / \n\,\OmegaS)
\;\cong\;  (\Lambda^2_S \OmegaS) /(\n\,\OmegaS \wedge \OmegaS)
\;\cong\; \langle 0\rangle
$$
since $\Lambda^2_S \OmegaS = \Omega^2_{S/K} = \langle 0\rangle $.
Consequently, the $S/\n$-vector space dimension of~$\n/\n^2$
is one, and hence~$\n$ is unigenerated.

The equivalence of (b) and (c) follows from Proposition~\ref{prop:HigherKDM}.
\end{proof}

Based on this characterization, we now derive explicit formulas
for the Hilbert polynomials of $\OmegaR$ and $\Omega^2_{R/K}$ in the case
of a weakly curvilinear scheme~$\X$. The following proposition provides
the key ingredients.

\begin{proposition} 
Let $(T,\n)$ be a 0-dimensional local affine $K$-algebra which is weakly
curvilinear. Write $\n = \langle a\rangle$ with $a\in T$, and
let~$\nu$ be the index of nilpotency of~$a$, i.e., let $\nu = \min \{i\ge 1 \mid a^i = 0\}$.
Furthermore, let $L=T/\n$ be the residue field of~$T$ and $\kappa = \dim_K(L)$.
Then the following claims hold.
\begin{enumerate}
\item[(a)] We have $\dim_K(T) = \dim_K (\gr_\n(T))_{\mathstrut} = \nu\, \kappa$.

\item[(b)] We have $\dim_K ( \Omega^1_{\gr_\n(T)/K} ) = 
\begin{cases}  \nu\,\kappa & \textrm{\ if\ }\charac(K) \mid \nu \\   
(\nu - 1)\, \kappa & \textrm{\ if\ } \charac(K) \nmid \nu
\end{cases}$.

\item[(c)] The {\bf Noether different} $\theta_N(T/K)$ of the algebra~$T/K$ satisfies
$$
\theta_N(T/K) \;=\; \Ann_T(\Omega^1_{T/K}) \;=\; 
\begin{cases}
\langle 0\rangle  & \textrm{\ if\ }\charac(K) \mid \nu \\  
\langle a^{\nu -1} \rangle & \textrm{\ if\ } \charac(K) \nmid \nu
\end{cases}
$$

\item[(d)] We have $\dim_K ( \Omega^1_{T/K} ) = \dim_K ( \Omega^1_{\gr_\n(T)/K} )$.

\end{enumerate}
\end{proposition}

\begin{proof}
To show~(a), it suffices to note that $\langle 0\rangle = \n^\nu \subseteq
\n^{\nu - 1} \subseteq \cdots \subseteq \n \subseteq T$ yields
$$
\dim_K(T) \;=\;  \dim_K(T/\n) + \dim_K (\n/\n^2) + \cdots + \dim_K(\n^{\nu - 1}/\n^\nu)
\;=\; \dim_K(\gr_\n(T)).
$$

Next we prove~(b). Letting $\bar{a} = a + \n^2$, we have $\gr_\n(T) = L[\bar{a}]
\cong L[z] / \langle z^\nu \rangle$, and therefore $\Omega^1_{\gr_\n(T)/K}
\cong L[z] dz / (L[z]\cdot z^\nu \, dz  + L[z]\cdot \nu \, z^{\nu - 1}\, dz)$.
So, if we set $\nu' = \nu$ if $\charac(K) \mid \nu$ and $\nu' = \nu - 1$ 
if $\charac(K) \nmid \nu$
then we get $\Omega^1_{\gr_\n(T)/K} \cong L[z]/ \langle z^{\nu'}\rangle \cong L^{\nu'}$,
and the claim follows.

The first equality in~(c) follows from~\cite{Kun1986}, Proposition 10.18.
The remaining claims were shown in~\cite{SS1974}, \S 18, Beispiel~3. For the convenience of the reader, 
we recall the main steps of the proof in our notation. Let $b_1,\dots,b_\kappa \in T$
be elements such that their residue classes $\bar{b}_1,\dots,\bar{b}_\kappa \in L$ form
a $K$-basis of~$L$. Without loss of generality we may assume that the canonical trace map
$\Tr_{L/K}:\, L \longrightarrow K$ satisfies $\Tr_{L/K}(\bar{b}_1)=1$ and
$\Tr_{L/K}(\bar{b}_i) = 0$ for $i=2,\dots,\kappa$.
Then the elements in $B= \{ a^i \, b_j \mid i\in \{0,\dots,\nu-1\},\, j \in \{1,\dots,\kappa\}\}$
form a $K$-basis of~$T$. Since we have $\n \cdot \langle a^i\rangle = \langle a^{i+1}\rangle$
for all $i\ge 0$, the canonical trace map $\Tr_{T/K}:\, T \longrightarrow K$ satisfies
$\Tr_{T/K}(c) = \nu \cdot \Tr_{L/K}(\bar{c})$ for every $c\in T$.

In the case $\charac(K) \mid \nu$, it follows that $\Tr_{T/K}=0$. As~$T$ is a Gorenstein
ring with a trace map~$\sigma$ and $\Tr_{T/K} =  g \cdot \sigma$ for a generating element
$g\in T$ of~$\theta_N(T/K)$ by~\cite{Kun1986}, Corollary F.12, it follows that 
$\theta_N(T/K) = \langle 0\rangle$, as claimed.

Now assume that $\charac(K) \nmid \nu$. Let $\sigma:\; T \longrightarrow K$ be the projection 
to~$a^{\nu-1}\, b_1$ along the rest of the basis~$B$. It is straightforward to check that~$\sigma$
is a trace map of~$T/K$ and that $\Tr_{T/K} = \nu\, a^{\nu-1}\, \sigma$. Therefore the element
$a^{\nu - 1}$ generates $\theta_N(T/K)$, as claimed.

It remains to prove~(d). From $\Omega^1_{T/K} = T\, da$ we get $\dim_K(\Omega^1_{T/K}) = 
\dim_K(T) - \dim_K(\theta_N(T/K))$. Thus the claim is a consequence of~(b) and~(c).
\end{proof}

As a consequence of this proposition and Proposition~\ref{prop:CharCurvilinearKD},
we have the following formulas for the Hilbert polynomials of $\OmegaRm$ 
when~$\X$ is weakly curvilinear.

\begin{corollary}\label{cor:ExplicitCurvilinear}
Let $\X$ be a 0-dimensional weakly curvilinear scheme in~$\mathbb{P}^n$,
let $S=\mathcal{O}_1 \times \cdots \times \mathcal{O}_t$ be the
decomposition of its affine coordinate ring into local rings,
let~$L_i$ be the residue class field of~$\mathcal{O}_i$, let $\kappa_i = 
\dim_K(L_i)$, and let $\nu_i = \dim_K(\mathcal{O}_i)/\kappa_i$ for $i=1,\dots,t$.
\begin{enumerate}
\item[(a)] We have $\HP(\OmegaR) = 2 \deg(\X) -
\textstyle{\sum\limits_{\charac(K)\nmid \,\nu_i}} \kappa_i$.

\item[(b)] We have $\HP(\Omega^2_{R/K}) = \deg(\X) 
- \textstyle{\sum\limits_{\charac(K)\nmid \,\nu_i}}  \kappa_i $.

\item[(c)] For $m>2$, we have $\HP(\OmegaRm) = 0$.
\end{enumerate}
\end{corollary}

\bigskip\bigbreak
%%%%%%%%%%%%%%%%%%%%%%%%%%%%%%%%%%%%%%%%%%%%%%%%%%%%%%%%%%%%
%
%  Section 7: K\"ahler Differentials of Fat Point Schemes
%
%%%%%%%%%%%%%%%%%%%%%%%%%%%%%%%%%%%%%%%%%%%%%%%%%%%%%%%%%%%%

\section{K\"ahler Differentials of Fat Point Schemes}%
\label{sec7:FatPointZD}

For fat point schemes, the above results about the Hilbert functions
and Hilbert polynomials of the modules of K\"ahler differentials 
can be made even more explicit.
Let $K$ be a perfect field and $P=K[X_0, \dots, X_n]$. Let $t\ge 1$, 
let $p_1,\dots,p_t$ be distinct  $K$-rational points in $\mathbb{P}^n$, 
and let $I_{p_i} \subseteq P$ be the homogeneous vanishing ideal of~$p_i$
for $i=1,\dots,t$. Given positive integers $m_1,\dots,m_t$,
recall that the 0-dimensional scheme~$\X$ defined by the saturated homogeneous ideal
$$
I_\X \;=\; I_{p_1}^{m_1} \;\cap\; \cdots \;\cap\; I_{p_t}^{m_t}
$$ 
is called the {\bf fat point scheme} with support $\Supp(\X) = \{p_1,\dots,p_t\}$
and {\bf multiplicities} $m_1,\dots,m_t$. Frequently, the scheme~$\X$ is written
as $\X = m_1 p_1 + \cdots + m_t p_t$. Let us collect some initial observations
about this setting.

\begin{remark}\label{rem:FatPointProps}
For a fat point scheme $\X = m_1 p_1 + \cdots + m_t p_t$ in $\mathbb{P}^n$, 
we know that:
\begin{enumerate}
\item[(a)] The degree of $\X$ is given by $\deg(\X)= \sum_{i=1}^t  \tbinom{n+m_i-1}{n}$.
	
\item[(b)] The scheme $\X$ is a reduced scheme if and only if
$\HP(\OmegaR)=t$ and $\HP(\OmegaRm)=0$ for $m\ge 2$.
\end{enumerate}
\end{remark}

It is natural to ask what the Hilbert polynomial of
$\Omega^m_{R/K}$ is for $1\le m\le n+1$. More precisely, does the Hilbert polynomial of
$\Omega^m_{R/K}$ depend only on $m$, $n$, and the multiplicties  $m_1,\dots,m_t$?

In view of Proposition~\ref{prop:HigherKDM}, in order to determine
the Hilbert polynomial of $\Omega^m_{R/K}$ for a fat point scheme~$\X$,
it suffices to work out the dimension of the $K$-vector spaces 
$\Omega^m_{S/K}$ for $m=1,\dots,n$. According to \cite{Kun1986}, Proposition~4.7, we have 
$$
\Omega^m_{S/K} = \Omega^{m}_{\mathcal{O}_1/K}\times\cdots\times
\Omega^{m}_{\mathcal{O}_t/K}
$$
and so $\dim_K(\Omega^m_{S/K}) = 
\sum_{i=1}^t\dim_K(\Omega^{m}_{\mathcal{O}_i/K})$.
This leads us to compute $\dim_K(\Omega^{m}_{\mathcal{O}_i/K})$
at a fat point of $\X$. 

Using a homogeneous linear change of coordinates, we assume that $\X=kp$ 
with $p=(1:0: ... :0)$ and $k\ge 1$.
Then we have $I_p=\langle X_1,\dots,X_n\rangle \subseteq P$.
Letting $A:=K[X_1,\dots,X_n]$ and $\q := \langle X_1,\dots,X_n \rangle$,
the local ring of~$\X$ at~$p$ is $S = A/\q^k$.

In the reduced case $k=1$, we have $\Omega^m_{S/K}=\langle 0\rangle$
for all $m\ge 2$. Thus we consider the case $k\ge 2$
now. We equip~$A$ with the standard grading and note that $\q^k$ is a homogeneous ideal.
Therefore~$S$ is a graded 0-dimensional affine $K$-algebra. Its Hilbert function 
is given by
$$
\HF_S(i) \;=\; \begin{cases}
\tbinom{n-1+i}{n-1} & \hbox{\rm for\ } 0\le i\le k-1 \\ 
\quad 0 & \hbox{\rm for\ } i\ge k
\end{cases}
$$ 
and thus $\dim_K(S)=\sum_{i\ge 0} \HF_S(i) = \binom{n+k-1}{n}$.

Our next task is to describe the Hilbert functions of the graded $S$-modules
$\Omega^m_{S/K}$ for all $m\ge 1$ explicitly. Notice that $\Omega^m_{S/K}=0$ for $m>n$.

\begin{proposition}\label{prop:K-HFOmegaSm}
In the above setting, let $1\le m\le n$. Then the Hilbert function of~$\Omega^m_{S/K}$ 
is given by 
$$
\HF_{\Omega^{m}_{S/K}}(i) \;=\; 
\begin{cases}
	\binom{n}{m}\binom{n+i-m-1}{n-1} & \mbox{if}\ i < m+k-1\\
	\tbinom{n}{m}\tbinom{n+k-2}{n-1}-\delta & \mbox{if}\ i = m+k-1\\
	0 & \mbox{if}\ i > m+k-1
\end{cases}
$$
where $\delta = \dim_K(d\q^k \wedge \Omega^{m-1}_{A/K})_{m+k-1}$.
\end{proposition}

\begin{proof}
By \cite{Kun1986}, Proposition~4.12, we have a homogeneous short exact sequence of graded
$A$-modules 
$$
0 \;\longrightarrow\;  (\q^k \; \Omega^{m}_{A/K} + d\q^k \wedge \Omega^{m-1}_{A/K})
\;\longrightarrow\; \Omega^{m}_{A/K}  \;\longrightarrow\;  \Omega^{m}_{S/K} \;\longrightarrow\; 0
$$
where $\Omega^{m-1}_{A/K}= A$ if $m=1$.
Here $\Omega^m_{A/K} = \bigoplus_{1\le i_1<\cdots<i_m\le n} 
A dX_{i_1} \wedge\cdots\wedge dX_{i_m}$ is a graded free $A$-module of 
rank $\binom{n}{m}$ with basis elements of degree~$m$. Moreover, we have  
$\q^k = \langle\, t_1,\dots, t_N \,\rangle$ where $\mathbb{T}^n_k =\{t_1,\dots,t_N\}$
and $N=\binom{n+k-1}{n-1}$. By Proposition~\ref{prop:PresentationOmegaRm}, we have
$$
d\q^k \;=\; \langle\, dt_1,\dots,dt_N \,\rangle + \langle\, t_idX_j \mid 
1\le i\le N;\; 1\le j\le n \,\rangle
$$
In particular, the graded $A$-module 
$\q^k \, \Omega^m_{A/K} + d\q^k \wedge \Omega^{m-1}_{A/K}$ is generated in degrees 
$m+k-1$ and $m+k$. So, for $i<m+k-1$, we have 
$$
\HF_{\Omega^{m}_{S/K}}(i) \;=\; \HF_{\Omega^{m}_{A/K}}(i) \;=\; 
\tbinom{n}{m}\tbinom{n+i-m-1}{n-1}
$$
Furthermore, for $i\ge m+k$, we have $A_{i-m}=(\q^k)_{i-m}$ and
$$
(\Omega^{m}_{A/K})_i \;=\; \tsum_{1\le i_1<\cdots<i_m\le n} 
A_{i-m} dX_{i_1}\wedge\cdots\wedge dX_{i_m}  
\;=\; (\q^k\Omega^{m}_{A/K} + d\q^k \wedge \Omega^{m-1}_{A/K})_i
$$
Consequently, we get $\HF_{\Omega^{m}_{S/K}}(i)=0$.
In the case $i=m+k-1$, we have 
\begin{align*}
\HF_{\Omega^{m}_{S/K}}(m+k-1) &\;=\;  \HF_{\Omega^{m}_{A/K}}(m+k-1)
- \dim_K(d\q^k \wedge \Omega^{m-1}_{A/K})_{m+k-1}\\
&\;=\; \tbinom{n}{m}\tbinom{n+k-2}{n-1}- \dim_K(d\q^k \wedge \Omega^{m-1}_{A/K})_{m+k-1}
\end{align*}
because $(\q^k \Omega^m_{A/K})_{m+k-1} = 0$.
\end{proof}

To describe the Hilbert function of $\Omega^{m}_{S/K}$ completely, 
it remains to compute $\delta = \dim_K(d\q^k \wedge \Omega^{m-1}_{A/K})_{m+k-1}$.
For that, we simplify the notation as follows.

\begin{remark}\label{rem:maprho}
For every $m\ge 1$, the map $\rho:\; \Omega^m_{A/K} \longrightarrow 
A(-m)^{\binom{n}{m}}$ given by $\rho(f dX_{i_1} \wedge \cdots \wedge dX_{i_m}) = 
f e_{i_1,\dots,i_m}$ for $f\in A$ and $1\le i_1 < \cdots < i_m \le n$ is an isomorphism
of graded $A$-modules. Note that the degrees of the standard basis vectors satisfy
$\deg(e_{i_1,\dots,i_m}) = m$ here.

To compute the desired number~$\delta$, we
have to consider the graded $A$-submodule $U = \rho(d\q^k \wedge \Omega^{m-1}_{A/K})$
of $A(-m)^{\binom{n}{m}}$. By Proposition~\ref{prop:K-HFOmegaSm}, we have
$U_i = \{0\}$ for $i< m+k-1$ and $U_i = (A(-m)^{\binom{n}{m}})_i$ for $i> m+k-1$.
Therefore the only interesting homogeneous component of~$U$ is $U_{m+k-1}$,
and its $K$-dimension is the number~$\delta$ we are looking for.
The vector space $U_{m+k-1}$ will be called the {\bf initial defining vector space}
for $\Omega^1_{S/K}$.
\end{remark}

In order to describe a $K$-basis of $U_{m+k-1}$, we use the module term ordering
$\sigma = {\tt DegRevLex-Pos}$ on $A(-m)^{\binom{n}{m}}$ (see \cite{KR2000},
Definition 1.4.16) and determine a system of generators of~$U_{m+k-1}$ having 
distinct leading terms. Here the standard basis vectors $e_{i_1,\dots,i_m}$ are ordered by 
using the lexicographic ordering on the tuples $(i_1,\dots,i_m)$.
The next three propositions set the stage for our main result.

\begin{proposition}\label{prop:LTgradti}
In the setting of the preceding remark, assume that $\charac(K)=0$
or $\charac(K)>k$.
Let $\{t_1,\dots,t_N\}$ be the set of terms of degree~$m+k-1$ in~$A$,
where $N = \binom{n+k-1}{n-1}$, and where the terms are ordered decreasingly with 
respect to {\tt DegRevLex}.
\begin{enumerate}
\item[(a)] For $i=1,\dots,N$, we have $\rho(dt_i) = \tfrac{\partial t_i}{\partial X_1}e_1 
+ \cdots +\tfrac{\partial t_i}{\partial X_n}e_n$. We denote the vector 
on the right-hand side by $\grad(t_i)$.

\item[(b)] Let $i\in \{1,\dots,N\}$ and write 
$t_i=X_1^{\alpha_1}\cdots X_n^{\alpha_n} \in \mathbb{T}^n_k$, where $\alpha_j \ge 0$.
Then we have   
$$
\LT_\sigma(\grad(t_i)) \;=\; 
\begin{cases}
   X_1^{\alpha_1}X_2^{\alpha_2}\cdots X_n^{\alpha_n-1} e_n & \mbox{if}\quad \alpha_n > 0\\
   X_1^{\alpha_1}X_2^{\alpha_2}\cdots X_{n-1}^{\alpha_{n-1}-1} e_{n-1} & 
        \mbox{if}\quad \alpha_n= 0,\alpha_{n-1}>0\\
   \vdots & \\
   X_1^{\alpha_1-1}e_{1} & \mbox{if}\quad \alpha_n=\cdots=\alpha_2= 0,\alpha_{1}=k
\end{cases}
$$

\item[(c)] For $i\in\{1,\dots,N\}$, we can write $\LT_\sigma(\grad(t_i))= 
X_1^{\beta_1}X_2^{\beta_2}\cdots X_\ell^{\beta_\ell}\, e_\ell$ with 
$\ell \in\{1,\dots,n\}$ and $\beta_j \ge 0$.
Then we have $t_i= X_1^{\beta_1} X_2^{\beta_2} \cdots X_\ell^{\beta_\ell+1}$.

Consequently, $\LT_\sigma(\grad(t_i))$ and~$t_i$ uniquely determine each other
and all leading terms $\LT_\sigma(\grad(t_i))$ are pairwise distinct.

\item[(d)] In the case $m=1$, the set $B_1 = \{\LT_\sigma(\grad(t_1)), \dots,
\LT_\sigma(\grad(t_N)) \}$ has $N = \binom{n+k-1}{n-1}$ distinct elements.

\end{enumerate}
\end{proposition}

\begin{proof}
Claim (a) follows from $dt_i = \tfrac{\partial t_i}{\partial X_1} X_1 
+ \cdots + \tfrac{\partial t_i}{\partial X_n} X_n$ and claim~(b)
is a consequence of 
$\grad(t_i) = \sum_{j=1}^n \alpha_j X_1^{\alpha_1}\cdots X_j^{\alpha_j-1}
\cdots X_n^{\alpha_n} e_j$. Claim~(c) follows from~(b) and claim~(d) 
follows from~(c).
\end{proof}

In the case $m\ge 2$, the determination of the leading terms in $U_{m+k-1}$
is slightly more involved.

\begin{proposition}{\bf (Generators of the Initial Defining Vector 
Space)}\label{prop:GensviJ}\\
In the setting described above, let $m\ge 2$, and assume
that $\charac(K)=0$ or $\charac(K)>k$. For $i\in \{1,\dots,N\}$ and 
$1\le j_1 < \cdots < j_{m-1} \le n$, we let $J=\{j_1,\dots,j_{m-1}\}$, and
we define 
$$
v_{i,J} = \rho ( d t_i \wedge dX_{j_1} \wedge \cdots \wedge dX_{j_{m-1}} )
\in U_{m+k-1}.
$$
\begin{enumerate}
\item[(a)]  The set~$G$ of all vectors $v_{i,J}$ generates the
$K$-vector space $U_{m+k-1}$.

\item[(b)] For a set $J= \{j_1,\dots,j_{m-1}\}$ as above and $\ell\notin J$, 
choose $\nu_\ell\in \{0,\dots,m-1\}$ such that 
$1\le j_1<\cdots<j_{\nu_\ell}< \ell <j_{\nu_\ell+1} < \cdots < j_{m-1}\le n$, and let
$\bar{e}_{\ell,J} = (-1)^{\nu_\ell}\, e_{j_1,\dots,j_{\nu_\ell}, \ell, j_{\nu_\ell+1},\dots,j_{m-1}}$. 
Then we have
$$
v_{i,J} = \tsum_{\ell\notin J} \tfrac{\partial t_i}{\partial X_\ell}\, \bar{e}_{\ell,J}.
$$
\end{enumerate}
\end{proposition}

\begin{proof}
Claim~(a) follows from the fact that the elements $dt_i \wedge dX_{j_1} \wedge
\cdots \wedge dX_{j_{m-1}}$ generate the vector space
$(\q^k \Omega^m_{A/K} + d \q^k \wedge \Omega^{m-1}_{A/K})_{m+k-1}$
since $(\q^k \Omega^m_{A/K})_{m+k-1} = 0$.

To show~(b), we calculate 
\begin{gather*}
dt_i \wedge dX_{j_1} \wedge \cdots \wedge dX_{j_{m-1}} \;=\; 
\left( \tsum_{\ell=1}^n \tfrac{\partial t_i}{\partial X_\ell} dX_\ell \right)
\wedge dX_{j_1} \wedge \cdots \wedge dX_{j_{m-1}}\hfill \\
\hfill \;=\; \tsum_{\ell\notin J} \tfrac{\partial t_i}{\partial X_\ell}(-1)^{\nu_\ell}
dX_{j_1} \wedge \cdots \wedge dX_{j_{\nu_\ell}} \wedge dX_\ell \wedge dX_{j_{\nu_\ell+1}}
\wedge \cdots \wedge dX_{j_m}
\end{gather*}
and apply the map~$\rho$.
\end{proof}

After describing a system of generators of the initial defining vector space,
we now turn to finding its leading term vector space with respect to~$\sigma$.
Then we are able to deduce its dimension via $\dim_K(U_{m+k-1}) = 
\dim_K ( \LT_\sigma(U_{m+k-1}) )$.

\begin{proposition}{\bf (Leading Terms of the Initial Defining Vector 
Space)}\label{prop:LTviJ}\\
In the setting described above, assume that $\charac(K)=0$ or
$\charac(K) > k$.
\begin{enumerate}
\item[(a)] Let $i\in\{1,\dots,N\}$ and $J=\{j_1, \dots,j_{m-1} \}$,
where $1\le j_1 < \cdots < j_{m-1} \le n$. Then we have
$$
\LT_\sigma(v_{i,J})  \;=\;  \tfrac{t_i}{X_\ell} \, e_{\ell,J}
$$
where $\ell \in \{1,\dots,n\}$ is the largest index such that
$\ell \notin J$ and $X_\ell \mid t_i$, and where $e_{\ell,J}
= e_{j_1,\dots,j_{\nu},\ell,j_{\nu+1},\dots,j_{m-1}}$.
If no such~$\ell$ exists, we have $v_{i,J}=0$.

\item[(b)] Let $i,i' \in \{1,\dots,N\}$ be such that $i\le i'$,
and assume that $J = \{j_1, \dots, j_{m-1}\}$ and $J' = \{j'_1,\dots, j'_{m-1}\}$
satisfy $1\le j_1 < \cdots < j_{m-1} \le n$ as well as
$1\le j'_1 < \cdots < j'_{m-1} \le n$. Suppose that we have $v_{i,J} \ne v_{i',J'}$
and $\LT_\sigma(v_{i,J}) = \LT_\sigma(v_{i',J'})$. Then the following
properties hold.

\begin{itemize}
\item[(1)] $i < i'$

\item[(2)] There exists indices $\ell\in \{1,\dots,n\}\setminus J$
and $\ell' \in \{1,\dots,n\} \setminus J'$
such that $t_i = X_1^{\alpha_1} \cdots X_n^{\alpha_n}$ and
$t_{i'} = X_1^{\alpha'_1} \cdots X_n^{\alpha'_n}$ satisfy
$\alpha_\ell = \alpha'_\ell +1$ and $\alpha'_{\ell'} = \alpha_{\ell'}+1$.

\item[(3)] $\ell < \ell'$

\item[(4)] The set $\hat{J} = J \cap J'$ satisfies 
$\hat{J} = J \setminus \{\ell'\} = J' \setminus \{ \ell \}$.

\item[(5)] We have decompositions $t_i = \hat{t} \cdot \tilde{t} \cdot X_\ell^{\alpha_\ell} 
X_{\ell'}^{\alpha_{\ell'}}$ and $t_{i'} = \hat{t} \cdot \tilde{t} \cdot 
X_\ell^{\alpha'_\ell} X_{\ell'}^{\alpha'_{\ell'}}$ where $\hat{t} \in 
K[X_j \mid j \in \hat{J}]$ and where
$\tilde{t} = X_{\nu_1}^{\alpha_{\nu_1}} \cdots X_{\nu_r}^{\alpha_{\nu_r}}$
with $\nu_1,\dots,\nu_r \notin \hat{J} \cup \{\ell, \ell' \}$.
\end{itemize}

\item[(c)] In the setting of~(b), we have
\begin{gather*}
v_{i,J} \;=\; \alpha_\ell\, \LT_\sigma(v_{i,J}) + 
\tsum_{\kappa=1}^r \hat{t}\; X_\ell^{\alpha_\ell} X_{\ell'}^{\alpha_{\ell'}}
\tfrac{\partial\tilde{t}}{\partial X_{\nu_\kappa}} \,
e_{\nu_\kappa, J}
\hbox{\quad \it and}\\
v_{i',J'} \;=\; \alpha'_{\ell'}\, \LT_\sigma(v_{i',J'}) + 
\tsum_{\kappa=1}^r \hat{t}\; X_\ell^{\alpha_\ell - 1} X_{\ell'}^{\alpha_{\ell'} + 1}
\tfrac{\partial\tilde{t}}{\partial X_{\nu_\kappa}}  \,
e_{\nu_\kappa, J'}
\end{gather*}

\item[(d)] The leading terms of the elements of~$G$ generate $\LT_\sigma(U_{m+k-1})$.

\item[(e)] The set of distinct $\sigma$-leading terms of elements of~$G$ is
\begin{align*}
B_m 
&:= \{\, t\, e_{j_1,\dots,j_m} \mid t = X_1^{\alpha_1}\cdots X_\ell^{\alpha_\ell} 
\in \mathbb{T}^\ell_{k-1}, 1\le j_1 < \cdots < j_m\le n,\, \ell\le j_m, \alpha_\ell>0 \}\\
&= \{ \LT_\sigma(v_{i,J}) \mid i\in \{1,\dots,N\},\, J=\{j_1,\dots,j_{m-1}\},\,
1 \le j_1 < \cdots < j_{m-1} \le n \}.
\end{align*}

\item[(f)] We have $\# B_m = \tbinom{n}{m}\tbinom{n+k-2}{n-1} 
- \tbinom{m+k-2}{m}\tbinom{n+k-2}{n-m-1}$.

\end{enumerate}
\end{proposition}

\begin{proof}
First we prove~(a). If $\ell\in J$ for all indices $\ell\in\{1,\dots,n\}$ such that
$X_\ell \mid t_i$ then Proposition~\ref{prop:GensviJ}.b shows $v_{i,J}=0$. 
So, let us assume that the largest index~$\ell$ with $X_\ell \mid t_i$ and $\ell\notin J$ exists. 
In view of Proposition~\ref{prop:GensviJ}.b, we can write $\LM_\sigma(v_{i,J})
= \tfrac{\partial t_i}{\partial X_\ell}\, \bar{e}_{\ell,J}$, and we only have to note that
$\tfrac{\partial t_i}{\partial X_\ell}$ differs from $\tfrac{t_i}{X_\ell}$
by a unit of~$K$.

Next we show the properties claimed in~(b) one by one. To verify (1), we remark
that $i=i'$ implies that the power products in $\LT_\sigma(v_{i,J})$ and
$\LT_\sigma(v_{i',J'})$ are equal. But then also the positions of these 
two leading terms have to be equal, i.e., we have $J=J'$. Altogether, the
equality $(i,J) = (i',J')$ contradicts $v_{i,J} \ne v_{i',J'}$.

To check~(2), we use~(a) and find $\ell,\ell' \in \{1,\dots,n\}$ such that
the power product in $\LT_\sigma(v_{i,J}) = \LT_\sigma(v_{i',J'})$ is
$\tfrac{t_i}{X_\ell} = \tfrac{t_{i'}}{X_{\ell'}}$. As we just saw, the two power products
cannot be equal, whence $\ell\ne \ell'$. Now $X_{\ell'}\, t_i = X_\ell \, t_{i'}$
yields the desired equalities for the exponents of the indeterminates.

In order to see why~(3) holds, we have to recall that the terms $t_1,\dots,t_N$
were ordered decreasingly with respect to~$\tt DegRevLex$. Hence $i<i'$ yields
$t_i >_{\tt DegRevLex} t_{i'}$. Since the two terms differ only in the exponents
of~$X_\ell$ and $X_{\ell'}$, it follows that $\ell<\ell'$, as the exponent of~$t_{i'}$
is larger for the indeterminate $X_{\ell'}$.

Claim~(4) follows immediately from~(2). Thus it remains to verify~(5).
Let $\tilde{t}_i$ be the power product of indeterminates of~$t_i$
which are not in $\hat{J}\cup \{\ell, \ell'\}$, and define $\tilde{t}_{i'}$ analogously.
Recall that $\ell = \max \{ \lambda \notin J \mid \alpha_\lambda >0\}$.
Therefore all indices $\lambda\notin \hat{J}\cup \{\ell, \ell'\}$ of indeterminates 
dividing~$\tilde{t}_i$ satisfy $\lambda <\ell$. Similarly, we have
$\ell' = \max\{ \lambda\notin \hat{J}\cup \{\ell\} \mid \alpha_\lambda >0 \}$.
This maximum is larger than~$\ell$, and the maximality of~$\ell$ implies that no such
indices~$\lambda$ exist between~$\ell$ and~$\ell'$. Hence all indices 
$\lambda\notin \hat{J}\cup \{\ell, \ell'\}$ of indeterminates 
dividing~$\tilde{t}_{i'}$ satisfy $\lambda <\ell$, and we obtain $\tilde{t}:= \tilde{t}_i =
\tilde{t}_{i'}$. Since the exponents of~$t_i$ and $t_{i'}$ differ only at~$X_\ell$
and $X_{\ell'}$, combining all indeterminates with indices in~$\hat{J}$ into~$\hat{t}$
yields the claimed decompositions of~$t_i$ and~$t_{i'}$.

To show~(c), we have to consider the image of $dt_i \wedge e_{j_1} \wedge \cdots
\wedge e_{j_{m-1}}$ under~$\rho$. Using the decomposition of~$t_i$ given in (b.5),
and considering the fact that $e_{\lambda} \wedge e_{j_1} \wedge \cdots
\wedge e_{j_{m-1}} =0$ for $\lambda \in J = \hat{J} \cup \{\ell'\}$, only the
terms involving the partial derivatives $\tfrac{\partial t_i}{\partial X_\lambda}$
with $X_\lambda \mid \tilde{t}$ or $\lambda=\ell$ survive the wedge product.
In view of the formula for $\LT_\sigma(v_{i,J})$ shown in~(a), this yields the claim
for~$v_{i,J}$. For~$v_{i',J'}$, it follows analogously.

Next we prove~(d). If in a linear combination $\sum_i c_i v_{i,J}$ with $c_i\in K$
the largest leading terms do not cancel, the result is a vector whose leading term
is one of the $\LT_\sigma(v_{i,J})$, and the claim holds. Now suppose that
the largest term appearing in the linear combination cancels out.
We use induction on the largest term appearing in the result of the
linear combination and show that it can be reduced using~$G$. Then it follows 
that the leading terms of the elements of~$G$ generate $\LT_\sigma(U_{n+k-1})$.

To start the induction, we note that the smallest term
in degree $m+k-1$ is $X_n^{k-1} e_{n-m+1,\dots,n}$. It is the leading term
of $\rho(dX_n^k\, e_{n-m+1} \wedge \cdots \wedge e_{n-1})$ and thus a leading term
in $\LT_\sigma(G)$. Subtracting the appropriate multiple of the corresponding
element of~$G$ has to yield zero, because no smaller term exists.

For the induction step, we note that any linear combination of the vectors
$v_{i,J}$ whose largest term cancels out can be written as a linear combination of
fundamental syzygies of the form
$$
\alpha_{\ell'} v_{i,J} - \alpha_{\ell} v_{i',J'} \;=\;
\alpha_{\ell'} \tsum_{\kappa=1}^r \hat{t}\; X_\ell^{\alpha_\ell} X_{\ell'}^{\alpha_{\ell'}}
\tfrac{\partial\tilde{t}}{\partial X_{\nu_\kappa}} \, e_{\nu_\kappa, J}
- \alpha_{\ell} \tsum_{\kappa=1}^r \hat{t}\; X_\ell^{\alpha_\ell - 1} X_{\ell'}^{\alpha_{\ell'} + 1}
\tfrac{\partial\tilde{t}}{\partial X_{\nu_\kappa}} \, e_{\nu_\kappa, J'}
$$
Thus it suffices to consider a fundamental syzygy
whose support contains the largest term for which the claim has not yet been shown.
As $X_\ell^{\alpha_\ell} X_{\ell'}^{\alpha_{\ell'}}$ is larger than $X_\ell^{\alpha_\ell - 1}
X_{\ell'}^{\alpha_{\ell'}+1}$, it follows that the leading term of this syzygy is
$\hat{t}\, X_{\ell'}^{\alpha_{\ell'}} \tfrac{\tilde{t}}{X_{\nu_r}}\, X_\ell^{\alpha_\ell}
e_{\nu_r,J}$. This term is the leading term of $v_{i'',J''}$ where
$t_{i''} = t_i \cdot X_{\ell'} / X_{\nu_r}$ and $J'' = \{\nu_r\} \cup \hat{J}$, because
we have $t_{i''} = \hat{t}\, X_{\nu_1}^{\alpha_{\nu_1}} \cdots X_{\nu_{r-1}}^{\alpha_{\nu_{r-1}}}
\cdot X_{\nu_r}^{\alpha_{\nu_r} -1}\, X_\ell^{\alpha_\ell} \, X_{\ell'}^{\ell'+1}$
and $\{\nu_r, \ell'\} \cup \hat{J} = \{\nu_r\} \cup J$. Consequently, if we subtract the
appropriate multiple of $v_{i'',J''}$, we get zero or a smaller leading term, and the
claim is a consequence of the induction hypothesis.

For the proof of~(e) we first show that every element of~$B_m$ is 
indeed a leading term of an element of~$G$.
Given $t\, e_{j_1,\dots,j_m}\in B_m$ with $t\in \mathbb{T}^n_k$, we let
$\tilde{t} = t\, X_{j_m}$ and $J= \{j_1,\dots,j_{m-1}\}$.
Then the degree of~$\tilde{t}$ is~$k$ and the hypothesis about the characteristic of~$K$
implies that $\alpha_{j_m} + 1 \le \deg(\tilde{t}) = k$ is a unit in~$K$. 
Choosing $i\in \{1,\dots,N\}$ such that $t_i = \tilde{t}$, we get
from~(a) that $\LT_\sigma(v_{i,J}) = t\, e_{j_1,\dots,j_m}$.

Conversely, part~(a) shows that every leading term of an element of~$G$ is in~$B_m$.
We also observe that the elements of~$B_m$ are clearly pairwise distinct.

Finally, to prove~(f), we count $\# B_m$. Given a fixed index $\ell\in\{1,\dots,n\}$, 
in order to form a term $t = X_\ell \cdot \hat{t}$ with $\hat{t}\in\mathbb{T}^{\ell}_{k-2}$, 
we have $\binom{\ell+k-3}{\ell-1}$ choices.
From all $\binom{n}{m}$ choices for $(j_1,\dots,j_m)$ with $j_1<\cdots<j_m$ we have to
subtract all $\binom{\ell-1}{m}$ choices where $\{j_1,\dots,j_m\} \subseteq \{1,\dots,\ell-1\}$.
Altogether, we get 
\begin{align*}
\# B_m &= \tsum_{\ell=1}^n \tbinom{\ell+k-3}{\ell-1}
[\tbinom{n}{m}-\tbinom{\ell-1}{m}]
= \tbinom{n}{m}\tsum_{\ell=0}^{n-1} \tbinom{\ell+k-2}{\ell}-
\tsum_{i=0}^{n-m-1} \tbinom{m+i+k-2}{m+i}\tbinom{m+i}{m}\\
&=\tbinom{n}{m}\tbinom{n+k-2}{n-1} - \!\!\! \tsum_{i=0}^{n-m-1} 
\tbinom{m+i+k-2}{m+k-2}\tbinom{m+k-2}{m}
= \tbinom{n}{m}\tbinom{n+k-2}{n-1} - \tbinom{m+k-2}{m}\tbinom{n+k-2}{n-m-1}
\end{align*}
and the proof is complete.
\end{proof}

Notice that in the key step of the proof of part~(d) we can even
be more explicit, as the following remark shows.

\begin{remark}
Let $v_{i,J}$ and $v_{i',J'}$ be two vectors such that $\LT_\sigma(v_{i,J}) = 
\LT_\sigma(v_{i',J'})$, where $i,i' \in \{1,\dots,N\}$ and $J,J'$ satisfy the
conditions of part~(b) of the proposition. Then we can write the fundamental syzygy of $v_{i,J}$
and $v_{i',J'}$ in the form
$$
\alpha_{\ell'} v_{i,J} - \alpha_{\ell} v_{i',J'} \;=\;
\tsum_{\kappa=1}^r \alpha_{\nu_\kappa} v_{\mu_\kappa, J_\kappa}
$$
where $v_{\mu_\kappa,J_\kappa}$ is the image of 
$\grad(X_{\ell'}\, t_i / X_{\nu_\kappa}) \wedge e_{\nu_\kappa,\hat{J}}$
under~$\rho$. This can be shown by a straightforward, but lengthy calculation
and provides another proof of~(d).
\end{remark}

Now we are ready to prove the main result of this section which
allows us to determine the Hilbert function of the K\"ahler differential
modules of a fat point.

\begin{theorem}\label{thm:K-DimOmegaSm}
Let $k\ge 1$, let~$K$ be a field of characteristic $\charac(K)=0$ or $\charac(K)>k$,
let $\q=\langle X_1,\dots,X_n\rangle$ be the homogeneous maximal ideal of
$A=K[x_1,\dots,x_n]$, and let $S=A/\q^k$.
\begin{enumerate}
\item[(a)] For every $m\ge 1$, we have 
$$
\delta \;=\;  \dim_K (d\q^k \wedge \Omega^{m-1}_{A/K})_{m+k-1} \;=\; 
\tbinom{n}{m}\tbinom{n+k-2}{n-1} - \tbinom{m+k-2}{m}\tbinom{n+k-2}{n-m-1}.
$$

\item[(b)] For every $m\ge 1$, we have
$$
\dim_K(\Omega^{m}_{S/K}) \;=\; \tbinom{n}{m}\tbinom{n+k-2}{n} +
\tbinom{m+k-2}{m}\tbinom{n+k-2}{n-m-1}.
$$
\end{enumerate}
\end{theorem}

\begin{proof}
To prove~(a), we identify $(d\q^k \wedge \Omega^{m-1}_{A/K})_{m+k-1}$
with the $K$-vector subspace $U_{m+k-1}$ of $(A(-m)^{\binom{n}{m}})_{m+k-1}$
using the map~$\rho$ of Remark~\ref{rem:maprho}. 

In the case $m=1$, we apply Proposition~\ref{prop:LTgradti} and get
$$
\dim_K(d\q^k)_k \;=\; \dim_K \langle\, \LT_\sigma(\grad(t_1)),\dots,
\LT_\sigma(\grad(t_N)) \,\rangle_K 
\;=\;  \tbinom{n+k-1}{n-1}
$$
Thus the claim follows from 
$$\tbinom{n+k-1}{n-1} = \tbinom{n+k-2}{n-1} + \tbinom{n+k-2}{n-2} =
n \tbinom{n+k-2}{n-1} - k \tbinom{n+k-2}{n-2} + \tbinom{n+k-2}{n-2}.
$$

Now let $m\ge 2$. In view of Proposition~\ref{prop:LTviJ}, it suffices 
to note that the elements of~$B_m$ form a $K$-basis of
$\LT_\sigma(U_{m+k-1})$ and to apply the well-known equality
$\dim_K(U_{m+k-1}) = \dim_K( \LT_\sigma(U_{m+k-1} ))$ 
(see for instance~\cite{KR2005}, Thm.~5.1.18).

To prove~(b) in the case $m=1$, we use Proposition~\ref{prop:K-HFOmegaSm} and get
$$
\dim_K(\Omega^1_{S/K})_k \;=\; n\tbinom{n+k-2}{n-1} -\tbinom{n+k-1}{n-1}
= (k-1)\tbinom{n+k-2}{n-2}.
$$
Therefore we have
\begin{align*}
\dim_K(\Omega^1_{S/K}) &\;=\; \tsum_{i=0}^{k-1}n\tbinom{n+i-2}{n-1} 
+ (k-1)\tbinom{n+k-2}{n-2} \\
&\;=\; n\tbinom{n+k-2}{n} + (k-1)\tbinom{n+k-2}{n-2}.
\end{align*}

Finally, we show~(b) in the case $m\ge 2$. 
By~(a) and Proposition~\ref{prop:K-HFOmegaSm}, we have
$$
\HF_{\Omega^m_{S/K}}(m+k-1) = \tbinom{n}{m}\tbinom{n+k-2}{n-1} -
\delta = \tbinom{m+k-2}{m}\tbinom{n+k-2}{n-m-1}.
$$
Consequently, we get
\begin{align*}
\dim_K(\Omega^m_{S/K}) &\;=\; \tsum_{i=0}^{m+k-2}\tbinom{n}{m}\tbinom{n+i-m-1}{n-1}
+ \tbinom{m+k-2}{m}\tbinom{n+k-2}{n-m-1} \\
&\;=\; \tbinom{n}{m}\tbinom{n+k-2}{n} + \tbinom{m+k-2}{m}\tbinom{n+k-2}{n-m-1}.
\end{align*}
Thus the proof of the theorem is complete.
\end{proof}

Under somewhat stricter assumptions on the characteristic of the ground field,
we can provide another proof of this theorem which is based on the Euler-Koszul complex
of the affine coordinate ring $S=R/\langle x_0-1\rangle = K[X_1,\dots,X_n]/\q^k$ of~$\X$,
where $\q = \langle X_1, \dots, X_n \rangle$. This complex is constructed as follows. 
In analogy to Definition~\ref{Euler-Form}, the algebra $S/K$ has the Euler form
$\epsilon_S: \Omega^1_{S/K} \rightarrow \q / \q^k$, which is given by 
$\epsilon_S(d_{S/K}f)=\deg(f)\cdot f$ for every homogeneous element $f\in S$.
After forming the Koszul complex of~$\epsilon_S$, we get the following result.

\begin{proposition}\label{prop:KoszulExactSeq}
In the above setting, let $\epsilon_S: \Omega^1_{S/K} \rightarrow \q/\q^k$
be the Euler form of~$S/K$.
\begin{enumerate}
\item[(a)] The sequence of $S$-linear maps
$$
0 \;\longrightarrow\; \Omega^n_{S/K} \;\TTo{\epsilon_S^{(n)}}\;
\cdots \; \TTo{\epsilon_S^{(3)}}\;  \Omega^2_{S/K}  \; \TTo{\epsilon_S^{(2)}}\;
\Omega^1_{S/K} \;\TTo{\epsilon_S}\; \q/\q^k \;\To\; 0 
$$
is a complex, called the Euler-Koszul complex of~$S/K$. Here $\epsilon_S^{(i)}$
is given by $\epsilon_S^{(i)}(w_1 \wedge \cdots \wedge w_i) = \sum_{j=1}^i
(-1)^{j-1} \epsilon_S(w_j)\, w_i \wedge \cdots \wedge \widehat{w_j} \wedge
\cdots \wedge w_i$ for $i\ge 2$ and $w_1,\dots,w_i \in \Omega^1_{S/K}$.

\item[(b)] Let $i\ge 1$, and assume that $\charac(K)=0$ or that $\charac(K)>0$ is
not a divisor of~$i$. Then the sequence of $K$-vector spaces
$$
0 \;\longrightarrow\; (\Omega^n_{S/K})_i    \;\TTo{\epsilon_S^{(n)}}\;
\cdots \; \TTo{\epsilon_S^{(3)}}\;  (\Omega^2_{S/K})_i  \; \TTo{\epsilon_S^{(2)}}\;
(\Omega^1_{S/K})_i \;\TTo{\epsilon_S}\; (\q/\q^k)_i \;\To\; 0
$$
is exact.
\end{enumerate}
\end{proposition}

\begin{proof}
Claim~(a) follows from the observation that this sequence is nothing but
the Koszul complex associated to the $S$-linear map~$\epsilon_S$.

To prove~(b) we have to show that $\Ker(\epsilon_S^{(m)}) \subseteq \Im(\epsilon_S^{(m+1)})$
for $m\ge 1$. In view of Proposition~\ref{prop:K-HFOmegaSm}, it is enough to check this
condition in the case $i\in \{1,\dots,m+k-1\}$.
It is straightforward verify that 
$(\epsilon_S^{(m+1)}\circ d + d\circ \epsilon_S^{(m)})(w)=iw$
for every $w\in (\Omega^m_{S/K})_i$. Given $w\in \Ker(\epsilon^{(m)})$,
we use the fact that $i\in K\setminus \{0\}$ and let
$w'=\frac{1}{i}dw \in (\Omega^{m+1}_{S/K})_i$. Then~$w'$ satisfies 
$\epsilon_S^{(m+1)}(w')= \frac{1}{i}(\epsilon^{(m+1)}\circ d 
+ d\circ \epsilon^{(m)})(w)= w$. Hence we have $w\in \Im(\epsilon^{m+1})$, and
the proof is complete.
\end{proof}

Based on this proposition, we can give an alternative proof of Theorem~\ref{thm:K-DimOmegaSm}.b
under a more stringent assumption about $\charac(K)$.

\begin{remark}\label{rem:AltProof}
In the setting of Theorem~\ref{thm:K-DimOmegaSm}, assume that
$\charac(K)=0$ or that $\charac(K) > m+k-1$. Then the proposition yields
\begin{align*}
\HF_{\Omega^m_{S/K}}(m+k-1) &\;=\; 
	\tsum_{\ell=m+1}^n (-1)^{\ell-m+1} \HF_{\Omega^\ell_{S/K}}(m+k-1)\\
	&\;=\; \tsum_{\ell=m+1}^n (-1)^{\ell-m+1}
	\tbinom{n}{\ell} \tbinom{n+m+k-\ell-2}{n-1}\\
	&\;=\; \tsum_{j=1}^{n-m} (-1)^{j+1}
	\tbinom{n}{m+j}\tbinom{n+k-j-2}{n-1} \;=\; \tbinom{m+k-2}{m}\tbinom{n+k-2}{n-m-1}
\end{align*}
where the last equality can be shown by induction. From this formula the
claim of Theorem~\ref{thm:K-DimOmegaSm}.b follows as in the proof of that theorem.
\end{remark}

The condition on the characteristic of~$K$ in Theorem~\ref{thm:K-DimOmegaSm} 
is necessary, as the following example shows.

\begin{example}
Consider the fat point scheme $\X = 2p$ 
in $\mathbb{P}^2$ over the field $K=\mathbb{F}_2$, where $p=(1:0:0)$. 
Then we have $k=2=\charac(K)$ and
\begin{align*}
\HF_S:\ 1\ 2\ 0\ 0\ \cdots,\quad 
\HF_{\Omega^1_{S/K}}:\ 0\ 2\ 3\ 0\ 0\cdots,\quad 
\HF_{\Omega^2_{S/K}}:\ 0\ 0\ 1\ 0\ 0\cdots
\end{align*}
This shows that $\dim_K(d\q^2)_2 = 1 \ne 3 =\tbinom{n}{1}\tbinom{n}{n-1} -
\tbinom{k-1}{1}\tbinom{n}{n-2}$ and 
$\dim_K(\Omega^1_{S/K})=5 \ne 3= \tbinom{n}{1}\tbinom{n}{n} +
\tbinom{k-1}{1}\tbinom{n}{n-2}$. In other words, the formulas 
given in Theorem~\ref{thm:K-DimOmegaSm} do not hold true for this case.

Moreover, we see that $\HF_S(2)-\HF_{\Omega^1_{S/K}}(2)+\HF_{\Omega^2_{S/K}}(2)\ne 0$.
Hence the sequence of $K$-vector spaces given 
in Proposition~\ref{prop:KoszulExactSeq} is not exact for $i=2$.
\end{example}

Our final result of this section provides an explicit formula
for the Hilbert polynomials of all K\"ahler differential modules
provided the characteristic of the base field is not too small.

\begin{theorem}\label{thm:HP-KDM}
Let $\X = m_1 p_1 + \cdots + m_t p_t$ be a fat point scheme in~$\mathbb{P}^n$,
let $m\in \{1,\dots,n+1\}$, and assume $\charac(K)=0$ or
$\charac(K) > \max\{m_1, \dots, m_t \}$.
Then we have 
$$
\HP(\Omega^m_{R/K}) = 
\begin{cases}
\tsum_{i=1}^t \; [\tbinom{n+m_i-1}{n} +(m_i-1)\tbinom{n+m_i-1}{n-1}] & \mbox{if}\quad m=1,\\
\tsum_{i=1}^t \; [\tbinom{n+1}{m}\tbinom{n+m_i-2}{n} +\delta_{i,m}] & \mbox{if}\quad 1<m\le n,\\
\tsum_{i=1}^t \tbinom{n+m_i-2}{n} & \mbox{if}\quad m=n+1. 
\end{cases}
$$
where $\delta_{i,m} = \tbinom{m+m_i-2}{m}\tbinom{n+m_i-2}{n-m-1}
+\tbinom{m+m_i-3}{m-1}\tbinom{n+m_i-2}{n-m}$.
\end{theorem}

\begin{proof}
By Proposition~\ref{prop:HigherKDM} and \cite{Kun1986}, Proposition~4.7, 
we have 
$$
\begin{aligned}
\HP(\Omega^{m}_{R/K}) &\;=\; \HP(\Omega^{m}_{\widetilde{R}/K})
= \dim_K(\Omega^{m}_{S/K}) + \dim_K(\Omega^{m-1}_{S/K})\\
&\;=\; \tsum_{i=1}^t (\dim_K(\Omega^{m}_{\mathcal{O}_i/K})
+\dim_K(\Omega^{m-1}_{\mathcal{O}_i/K}))
\end{aligned}
$$
For $m=1$, Theorem~\ref{thm:K-DimOmegaSm} yields
\begin{align*}
\dim_K(\mathcal{O}_i) + \dim_K(\Omega^{1}_{\mathcal{O}_i/K})
&\;=\; \tbinom{n+m_i-1}{n} + n\tbinom{n+m_i-2}{n} + (m_i-1)\tbinom{n+m_i-2}{n-2}\\
&\;=\; \tbinom{n+m_i-1}{n} +(m_i-1)\tbinom{n+m_i-1}{n-1}
\end{align*}
For $2 \le m \le n+1$, Theorem~\ref{thm:K-DimOmegaSm} shows 
$$
\begin{aligned}
\dim_K(\Omega^{m}_{\mathcal{O}_i/K})&+\dim_K(\Omega^{m-1}_{\mathcal{O}_i/K})\\
&\;=\; \tbinom{n+1}{m}\tbinom{n+m_i-2}{n} +\tbinom{m+m_i-2}{m}\tbinom{n+m_i-2}{n-m-1}
+\tbinom{m+m_i-3}{m-1}\tbinom{n+m_i-2}{n-m}
\end{aligned}
$$
In particular, we have 
$\dim_K(\Omega^{n+1}_{\mathcal{O}_i/K})+\dim_K(\Omega^{n}_{\mathcal{O}_i/K})
=\dim_K(\Omega^{n}_{\mathcal{O}_i/K})=\tbinom{n+m_i-2}{n}$.
Therefore we obtain the desired formula for $\HP(\Omega^m_{R/K})$.
\end{proof}

This theorem generalizes and unifies a number of special cases
which were considered previously. Let us mention one of them.

\begin{remark}
In the case $m=n+1$, Theorem~\ref{thm:HP-KDM} yields 
$$
\HP(\Omega^{n+1}_{R/K}) = \deg(\Y)= 
{\textstyle \sum\limits_{i=1}^t} \tbinom{m_i+n-2}{n}
$$
where $\Y$ is the fat point scheme 
$\Y = (m_1-1)\, p_1 + \cdots + (m_t-1)\, p_t$ in $\mathbb{P}^n$. 
This provides a positive answer to a question posed in~\cite{KLL2019}. 
See also~\cite{KLL2021} for a different proof of this result.
\end{remark}

\bigskip\bigbreak
%%%%%%%%%%%%%%%%%%%%%%%%%%%%%%%%%%%%%%%%%%%%%%%%%%%%%%%%%%%%
%
%  Section 8: K\"ahler Differentials of Uniform Schemes
%
%%%%%%%%%%%%%%%%%%%%%%%%%%%%%%%%%%%%%%%%%%%%%%%%%%%%%%%%%%%%

\section{K\"ahler Differentials of Uniform Schemes} \label{sec6:UniformZD}

In this section we examine how uniformity properties of 0-dimensional
schemes in~$\mathbb{P}^n$, in particular the Cayley-Bacharach property,
are reflected in the structure of their K\"ahler differential modules.
Thus we continue to let $K$ be a perfect field, let~$\X$ be a 0-dimensional 
subscheme of $\mathbb{P}^n$ having $K$-rational support,
and we let $R=P/I_\X$ be the homogeneous coordinate ring of~$\X$.
Recall that the Cayley-Bacharach property of~$\X$ is defined as follows.

\begin{definition}
Let $d\in\mathbb{N}$, and let $r_\X$ be the regularity index of~$\X$.
\begin{enumerate}
\item[(a)] The scheme~$\X$ is said to have 
the {\bf Cayley-Bacharach property} of degree~$d$ 
(in short, $\X$ has {\bf CBP($d$)}) 
if every hypersurface of degree~$d$ which contains a
subscheme~$\Y$ of~$\X$ of degree $\deg(\Y)= \deg(\X)-1$ 
automatically contains~$\X$.

\item[(b)] If $\X$ has the Cayley-Bacharach property of degree $r_\X-1$, 
then~$\X$ is called a {\bf Cayley-Bacharach scheme}.
\end{enumerate}
\end{definition}

Algebraically, the Cayley-Bacharach property can be characterized as follows.

\begin{remark}\label{rem:FirstCBP}
Let $\Y$ be a subscheme of~$\X$ of degree $\deg(\X)-1$. 
\begin{enumerate}
\item[(a)] The image of the vanishing ideal $I_\Y \subseteq P$ of~$\Y$ in $R=P/I_\X$ 
is denoted by $I_{\Y/\X} = I_\Y / I_\X$. Its initial degree
$\alpha_{\Y/\X} = \min\{ k\in \mathbb{N} \mid (I_{\Y/\X})_k \ne 0\}$
is called the {\bf separator degree} of~$\Y$ in~$\X$.
Then the Hilbert function of~$\Y$ satisfies
$$
\HF_\Y(i) \;=\;  \begin{cases}
\HF_\X(i) & \hbox{for}\quad i<\alpha_{\Y/\X},\\ 
\HF_\X(i)-1 & \hbox{for}\quad i\ge \alpha_{\Y/\X},
\end{cases}
$$
and we have $\alpha_{\Y/\X}\le r_\X$.

\item[(b)] A non-zero homogeneous element $f^\ast_\Y \in I_{\Y/\X}$
of degree $\alpha_{\Y/\X}$ is called a {\bf minimal separator} of~$\Y$ in~$\X$.
This case we have $(I_{\Y/\X})_{\alpha_{\Y/\X}+i} = K\, x_0^i\, f^\ast_\Y$
for every $i\ge 0$. 

\item[(c)] Using this terminology, the following conditions are equivalent
for every $d\ge 1$ (see also \cite{Kre1994}, Proposition~2.1):
\begin{itemize}
\item[(1)] $\X$ has ${\rm CBP}(d)$.

\item[(2)] Every subscheme $\Y$ of $\X$ of degree $\deg(\X)-1$
satisfies $\alpha_{\Y/\X} \ge d+1$.

\item[(3)] No element of $(I_{\Y/\X})_{r_\X} \setminus \{0\}$ 
is divisible by $x_0^{r_\X-d}$.

\end{itemize}
\end{enumerate}
\end{remark}

Notice that the number $r_{\X}-1$ is the largest degree $d\ge0$
such that $\X$ can have CBP($d$). Thus~$\X$ is a Cayley-Bacharach scheme
if and only if it has the CBP with respect to the largest possible degree.
The K\"ahler differential module of a subscheme~$\Y$ of~$\X$
of degree $\deg(\Y)=\deg(\X)-1$ can be described as follows.

\begin{proposition}\label{prop:OmegaSubsch} 
Let $\Y$ be a proper subscheme of~$\X$, 
and let $R_{\Y} = R/I_{\Y/\X}$ be the homogeneous coordinate ring of~$\Y$.
\begin{enumerate}
\item[(a)] We have $\Omega^1_{R_{\Y}/K} \cong \Omega^1_{R/K} / dI_{\Y/\X}$.
Moreover, for a system of non-zero homogeneous generators $\{g_1,\dots,g_r\}$ 
of~$I_{\Y/\X}$, we have
$$
dI_{\Y/\X} \;=\; \langle dg_i,\; g_i\,dx_j \mid i=1,\dots,r;\, j=0,\dots,n\rangle .
$$

\item[(b)] If $\deg(\Y)=\deg(\X)-1$ then we have
$dI_{\Y/\X} = \langle df_{\Y}^\ast,\, f_{\Y}^\ast\, dx_0 \rangle$,
where $f_{\Y}^\ast$ is a minimal separator of~$\Y$ in~$\X$.

\item[(c)] If $\charac(K)=0$ or $\charac(K) > r_{\X}$ then we have
$(dI_{\Y/\X})_{\alpha_{\Y/\X}} \ne 0$.
\end{enumerate}
\end{proposition}

\begin{proof} Claim~(a) is a consequence of 
Proposition~\ref{prop:PresentationOmegaRm}.d. Claim~(b) follows from
$(I_{\Y/\X})_k = K\, x_0^{k-\alpha_{\Y/\X}}\, f_{\Y}^\ast$ for
$k\ge \alpha_{\Y/\X}$.

It remains to prove~(c). Let $g\in I_{\Y/\X}$ be a non-zero homogeneous 
generator of least degree $\alpha_{\Y/\X}$. Since~$g$ is homogeneous, its
image $g\deh$ is non-zero in the affine coordinate ring $S=R/\langle x_0-1\rangle$
of~$\X$. By~(a), it suffices to show $dg\ne 0$.

Let $\Phi:\; \OmegaR \To \OmegaRtilde \cong
S[x_0] dx_0 \oplus K[x_0] \otimes_K \OmegaS $ 
be the canonical $K[x_0]$-linear map constructed
in Proposition~\ref{prop:CompareKahlerMod}.
The composition of~$\Phi$ with the projection to the first summand
$\Theta: \OmegaR \stackrel{\Phi}{\rightarrow} \OmegaRtilde
\rightarrow S[x_0] dx_0$ satisfies 
$\Theta(f dx_i) = f\deh x_i x_0^k dx_0$ for $i=1,\dots,n$
and $\Theta(f dx_0) = f\deh x_0^k dx_0$ for $f\in R_k$.
Using $dg = \frac{\partial g}{\partial x_0} dx_0+
\frac{\partial g}{\partial x_1} dx_1 +\cdots+
\frac{\partial g}{\partial x_n} dx_n$, we obtain
$$
\begin{aligned}
\Theta(dg) &= \Theta(\tfrac{\partial g}{\partial x_0} dx_0) +
\Theta(\tfrac{\partial g}{\partial x_1} dx_1)+ \cdots+
\Theta(\tfrac{\partial g}{\partial x_n} dx_n)\\
&= ((\tfrac{\partial g}{\partial x_0})\deh +
(\tfrac{\partial g}{\partial x_1})\deh x_1+ \cdots+
(\tfrac{\partial g}{\partial x_n})\deh x_n)\; x_0^{\alpha_{\Y/\X}-1}dx_0\\
&= (\tfrac{\partial g}{\partial x_0}x_0 +
\tfrac{\partial g}{\partial x_1} x_1+ \cdots+
\tfrac{\partial g}{\partial x_n} x_n)\deh\,  x_0^{\alpha_{\Y/\X}-1}dx_0\\
&=\alpha_{\Y/\X}\; g\deh \,x_0^{\alpha_{\Y/\X}-1}\,dx_0
\end{aligned}
$$
where the last equality follows from Euler's relation.
Since $\alpha_{\Y/\X}\ne 0$ in $K$ and $g\deh \ne 0$ in $S$,
we get $\Theta(dg)\ne 0$. Therefore we have $dg \ne 0$, as desired. 
\end{proof}

This proposition yields the following characterization of the Cayley-Bacharach 
property of~$\X$ using modules of K\"ahler differentials.

\begin{corollary}\label{cor:CBP-Omega}
Suppose that $\charac(K)=0$ or $\charac(K)>r_\X$.
Then the scheme~$\X$ has CBP($d$) if and only if
$\HF_{\OmegaR}(d) = \HF_{\Omega^1_{R_\Y/K}}(d)$
for every subscheme $\Y \subseteq \X$ of degree $\deg(\Y)=\deg(\X)-1$. 
\end{corollary}

\begin{proof}
First we assume that $\X$ has CBP($d$). 
Then $\alpha_{\Y/\X}\ge d+1$ by Remark~\ref{rem:FirstCBP}.c,
and hence the homogeneous generators of $dI_{\Y/\X}$ 
have degree $\ge d+1$. This implies 
$\HF_{\OmegaR}(d) = \HF_{\Omega^1_{R_\Y/K}}(d)$.

Conversely, assume that $\X$ does not have CBP($d$) and
that $\Y$ is a subscheme of~$\X$ with $\deg(\Y)=\deg(\X)-1$ and $\alpha_{\Y/\X}\le d$.
Then Proposition~\ref{prop:OmegaSubsch}.c yields 
$(dI_{\Y/\X})_{\alpha_{\Y/\X}}\ne 0$, and hence 
$\HF_{\OmegaR}(d) > \HF_{\Omega^1_{R_\Y/K}}(d)$.
\end{proof}

Recall our assumption that the scheme~$\X$ has $K$-rational support.
In this way we make sure that there exist subschemes~$\Y$ of~$\X$
of degree $\deg(\Y)=\deg(\X)-1$, because at each point in the support 
of~$\X$ the local ring contains a socle element that generates
an ideal which is a 1-dimensional $K$-vector space.

In Sections~2 and~3 of~\cite{KLR2019}, the notion of a Cayley-Bacharach 
scheme was generalized to arbitrary 0-dimensional schemes~$\X$ in $\mathbb{P}^n$.
The following example shows that Corollary~\ref{cor:CBP-Omega} fails to hold 
when the scheme~$\X$ does not have $K$-rational support.

\begin{example}
Consider the 0-dimensional scheme $\X$ in $\mathbb{P}^2$ over 
$\mathbb{Q}$ defined by 
$I_\X=\langle (X_1-X_0)^2, X_2^3+2X_0^2X_2+X_0^3\rangle$. 
Then $\deg(\X)=6$ and $\X$ does not have $K$-rational support. 

The scheme~$\X$ does not have subschemes of degree $\deg(\X)-1$.
Instead, it has a unique maximal subscheme~$\Y$ of degree~$3$,
namely the scheme defined by
$I_{\Y}= \langle X_1-X_0, X_2^3+2X_0^2X_2+X_0^3\rangle$. 
We compute $\HF_\X:\ 1\ 3\ 5\ 6\ 6\cdots$,
$r_\X=3$, and $\HF_\Y:\ 1\ 2\ 3\ 3\cdots$. 
Thus the scheme~$\X$ is a Cayley-Bacharach 
scheme in the sense of \cite{KLR2019}, Definition~3.10. 

However, further calculations show
that we have $\HF_{\OmegaR}: 0\ 3\ 8\ 12\ 12\ 10\ 9\ 9\cdots$
and  $\HF_{\Omega^1_{R_{\Y}/K}}: 0\ 2\ 4\ 5\ 4\ 3\ 3\cdots$.
In particular, we have $\HF_{\OmegaR}(2) \ne \HF_{\Omega^1_{R_{\Y}/K}}(2)$,
and hence Corollary~\ref{cor:CBP-Omega} does not hold.
\end{example}

The Cayley-Bacharach property can be interpreted as a weak uniformity
condition on the scheme~$\X$. It is a special case of the following more
general property.

\begin{definition}
Let $1\le i<\deg(\X)$ and $1\le j\le r_\X - 1$. We say that the scheme
$\X$ is \textbf{$(i,j)$-uniform}, if every subscheme $\Y\subseteq \X$ 
of degree $\deg(\Y)=\deg(\X)-i$ satisfies $\HF_\Y(j)=\HF_\X(j)$.
\end{definition}

In this terminology, the scheme $\X$ has CBP($d$) iff
it is $(1,d)$-uniform. Another well-known case is $(\deg(\X)-n-1,1)$-uniformity
which is also called {\it linearly general position}. If all possible
$(i,j)$-uniformities are satisfied, the scheme~$\X$ is commonly said to
be {\it in uniform position}.
Using Proposition~\ref{prop:OmegaSubsch}.c, we can characterize $(i,j)$-uniformity
in terms of Hilbert functions of K\"ahler differential modules as follows.

\begin{corollary}\label{cor:Uniformity-Omega}
Suppose that $\charac(K)=0$ or $\charac(K)>r_\X$.
Then the scheme $\X$ is $(i,j)$-uniform if and only if
$\HF_{\OmegaR}(j) = \HF_{\Omega^1_{R_\Y/K}}(j)$
for every subscheme $\Y \subseteq \X$ of degree $\deg(\Y)=\deg(\X)-i$. 
\end{corollary}

\begin{proof}
For a subscheme~$\Y$ of~$\X$ of degree $\deg(\Y)=\deg(\X)-i$,
we have the exact sequence
$$
0 \longrightarrow dI_{\Y/\X} 
\longrightarrow \OmegaR \longrightarrow \Omega^1_{R_\Y/K}
\longrightarrow 0.
$$
Clearly, the condition $\HF_\Y(j)=\HF_\X(j)$ is equivalent to 
$\alpha_{\Y/\X}>j$. By Proposition~\ref{prop:OmegaSubsch}.c, 
this is in turn equivalent to $(dI_{\Y/\X})_j=\langle 0\rangle$, 
and hence to $\HF_{\OmegaR}(j) = \HF_{\Omega^1_{R_\Y/K}}(j)$.
\end{proof}

%%%%%%%%%%%%%%%%%%%%%%%%%%%%%%%%%%%%%%%%%%%%%%%%%%%%
%
%  Bibliography
%
%%%%%%%%%%%%%%%%%%%%%%%%%%%%%%%%%%%%%%%%%%%%%%%%%%%%

\bigbreak

\end{document}